\documentclass[10pt,twoside]{article}
\usepackage{pifont}
\usepackage{bbm}
\usepackage{hyperref}
\usepackage{amsfonts}
\usepackage{amssymb}
\usepackage{amsmath}
\usepackage{amscd}
\usepackage{mathrsfs}
\input amssym.def
\input amssym.tex
\usepackage[all]{xy}
\def\bc{\begin{center}}
\def\ec{\end{center}}
\def\no{\noindent}

\def\D{\Delta}

\def\g{\gamma}

\def\o{\otimes}

\def\r{\rho}

\def\s{\sigma}

\def\v{\varepsilon}

\begin{document}
\begin{center}
{\Large\bf Hopf brace, braid equation and bicrossed coproduct$^{\ast}$}

\vspace{4mm}

{\bf Huihui Zheng$^1$, Fangshu Li$^1$, Tianshui Ma$^2$, Liangyun Zhang$^{1,\ast\ast}$}

~$^1$College of Science, Nanjing Agricultural University,
Nanjing 210095, China

$^2$College of Mathematics and Information Science, Henan Normal University,
Xinxiang 453007, China

\begin{figure}[b]
\rule[-2.5truemm]{5cm}{0.1truemm}\\[2mm]
\no $^{\ast}$This work is supported by Natural Science Foundation of China (11571173).

$^{\ast\ast}$Corresponding author: zlyun@njau.edu.cn
\end{figure}
\begin{minipage}{144.6mm}

\vskip1cm \footnotesize{\textbf{Abstract}:
In this paper, we mainly give some equivalent characterisations of Hopf braces, show that the category $\mathcal{CB}(A)$ of Hopf braces is equivalent to the category $\mathcal{C}(A)$ of bijective 1-cocycles, and prove that the category $\mathcal{CB}(A)$ of Hopf braces is also equivalent to the category  $\mathcal{M}(A)$ of Hopf matched pairs. Moreover, we construct many more Hopf braces on polynomial Hopf algebras, Long copaired Hopf algebras and Drinfel'd doubles of finite dimensional Hopf algebras, and give a sufficient and necessary condition for a given bicrossed coproduct $A\bowtie H$ to be a Hopf brace if $A$ or $H$ is a Hopf brace.
\vspace{3mm}}\vspace{3mm}\\
\textbf{2010 AMS Classification Number}: 16T05
\end{minipage}
\end{center}

\vspace{4mm}

\begin{center}
{\bf \S1\quad Introduction and Preliminaries}
\end{center}

Braces were introduced in \cite{W. Rump2} by Rump, which are a generalization of Jacobson radical rings, to understand the structure behind non-degenerate involutive set-theoretic solutions of Yang-Baxter equations. They provide a powerful algebraic framework to work with set-theoretic solutions and have also an advantage to discuss braided groups and sets imitating ring theory. Moreover, they have also connections with regular subgroups and orderable groups \cite{D. Bachiller3}, flat manifolds \cite{W. Rump2014}, Hopf-Galois extensions \cite{D. Bachiller2016}. And then they were generalized to skew braces by Guarnieri and Vendramin in \cite{L. Guarnieri}. Skew braces are used to study regular subgroups and Hopf-Galois extensions, bijective 1-cocycles, and triply factorized groups, see for example \cite{A. Smoktunowiczi1}. Through their connection with Yang-Baxter equation and group theory, braces have attracted a lot of attention and obtained a wide range of more influential results, for examples \cite{A. L. Agore, A. Smoktunowiczi, D. Bachiller, D. Bachiller2016, D. Bachiller2018, D. Bachiller1, D. Bachiller2, F. Ced, T. Gateva-Ivanova, W. Rump1, W. Rump2007, W. Rump2014}.

In \cite{I. Angiono}, the authors introduced the conception of Hopf braces. That is, a Hopf brace is a kind of special Hopf algebra with two different multiplications connected with antipode, which is of a new algebraic structure related to the Yang-Baxter equation, and is also a generalization of braces and skew braces.

It is well known that every Hopf algebra has an algebraic structure with a multiplication and a coalgebraic structure with a comultiplication. Thus, by the above discussion, we will introduce a new kind of Hopf algebra (also called a Hopf brace), with two different comultiplications connected with antipode.

The bicrossed product first emerged in group theory, which is constructed from a matched pair of groups, as a natural generalization of the semi-direct product.
Because of the connection between bicrossed product and factorization problem, the bicrossed product was considered and stuided in varied contexts, such as (co)algebras \cite{ A. Cap, S. Caenepeel}, Lie algebras \cite{S. Majid1, Y. Kosmann-Schwarzbach}, Hopf algebras \cite{S. Majid} and so on. The construction of bicrossed products in the context of Hopf algebras was introduced by Majid \cite{S. Majid}. We could  find that the Drinfel'd double $D(H)$ of a finite dimensional Hopf algebra $H$ was a special type of this bicrossed product in \cite{S. Majid}. As an important kind of associative algebras, it was closely related to smash product, L-R smash product, bismash product, crossed product and so on. Similarly, the bicrossed coproduct has the same effect as its dual form which is given by Hopf matched pair.

The main purpose of this paper is to construct Hopf braces with two comultiplications connected with antipode. The paper is organized as follows. In section 2, we give many examples of Hopf braces on polynomial Hopf algebras and Long copaired Hopf algebras, and prove that the category $\mathcal{CB}(A)$ of Hopf braces is equivalent to the category of bijective 1-cocycles $\mathcal{C}(A)$ given in Theorem 2.12. Moreover, we obtain a solution of the braid equation by a commutative Hopf brace.
 In section 3, we mainly study structures of commutative Hopf braces, build a correspondence between Hopf braces and Hopf matched pairs given in Proposition 3.3 and 3.4, and prove that the category $\mathcal{CB}(A)$ of Hopf braces is equivalent to the category $\mathcal{M}(A)$ of Hopf matched pairs $(A,A)$ given in Theorem 3.5.
 In section 4, we mainly give a sufficient and necessary condition for a given bicrossed coproduct $A\bowtie H$ to be a Hopf brace if $A$ or $H$ is a Hopf brace, and show that the dual of Drinfel'd double $D(H)$ of a finite dimensional cocommutative Hopf algebra $H$ is a Hopf brace.

 Throughout this paper, let $k$ be a fixed field, and our considered objects be all meant over $k$. And we freely use coalgebras, bialgebras, Hopf algebras, comodule algebras and comodule bialgebras terminology introduced in \cite{Sweedler, S. Montgomery} and \cite{Y. Chen}. For a coalgebra $C$, we write its comultiplication $\Delta(c)$ with $c_{1}\otimes c_{2}$, for any $c\in C$; for a left
$C$-comodule $M$, we denote its coaction by $\rho(m)=m_{(-1)}\otimes m_{(0)}$, for any $m\in M$; for a right
$C$-comodule $M$, we denote its coaction by $\rho(m)=m_{[0]}\otimes m_{[1]}$, for any $m\in M$, in which we omit the summation symbols for convenience.

\vspace{3mm}

\begin{center}
{\bf \S2\quad Hopf brace and its category}
\end{center}

\vspace{3mm}

In this section, we will introduce the conception of Hopf braces, give many examples of Hopf braces on polynomial Hopf algebras and Long copaired Hopf algebras, and mainly prove that the two category $\mathcal{CB}(A)$ of Hopf cobraces is equivalent to the full subcategory of the category $\mathcal{C}(A)$ of bijective 1-cocycles.

\vspace{2mm}

\textbf{Definition 2.1} \ \ Let $(H,m,1)$ be an algebra. A {\it Hopf brace} structure over $H$ consists of the following data:

(1) a Hopf algebra structure $(H, m, 1, \Delta, \varepsilon, S)$,

(2) a Hopf algebra structure $(H, m, 1, \Delta^\prime, \epsilon, T)$,

(3) satisfying the following compatibility:
$$
h_{1^\prime}\otimes h_{2^\prime1}\otimes h_{2^\prime2}=h_{11^\prime}S(h_2)h_{31^\prime}\otimes h_{12^\prime}\otimes h_{32^\prime}\eqno(2.1)
$$
for any $h\in H,$ where $\Delta(h)$ is denoted by $h_1\otimes h_2$ and $\Delta^\prime(h)$ denoted by $h_{1^\prime}\otimes h_{2^\prime}$.

\vspace{2mm}

\textbf{Remark 2.2}\ \ (1) In any given Hopf brace $(H, \Delta, \varepsilon; \Delta^\prime, \epsilon)$, we know that $\varepsilon=\epsilon$.

In fact, we applying $id\otimes \varepsilon\otimes \epsilon$ into the two sides of Eq.(2.1), we easily get
$$
h=h_{11^\prime}S(h_2)h_3\varepsilon(h_{12^\prime})=h_{1^\prime}\varepsilon(h_{2^\prime}).
$$

So, we obtain that $\epsilon(h)=\epsilon(h_{1^\prime}\varepsilon(h_{2^\prime}))=\epsilon(h_{1^\prime})\varepsilon(h_{2^\prime})=\varepsilon(h).$

(2) Let $(H, \Delta, \varepsilon)$ be a Hopf algebra. Then, we easily know that $(H, \Delta, \varepsilon; \Delta, \varepsilon)$ is a Hopf brace.

(3) Let $(H, \Delta, \varepsilon)$ be a Hopf algebra. If $H$ has not other Hopf braces unless itself as in (2), we say that it has not non-trivial Hopf braces.

 Let $U(L)$ be an universal enveloping algebra of a Lie algebra $L$. Then, $U(L)$ is a Hopf algebra, and $U(L)$ has not non-trivial Hopf braces since its comultiplication is unique.

(4) In what follows, we denote the Hopf brace in Definition 2.1 by $(H,\Delta,\Delta^\prime)$.

\vspace{2mm}

\textbf{Example 2.3}\ \ (1) Let $(H, \Delta, \varepsilon)$ be a Hopf algebra. Then, $(H,\Delta,\Delta^{coH})$ and $(H,\Delta^{coH},\Delta)$ are Hopf braces.

 (2) Let $A=k[g,g^{-1},x]$ be a Hopf algebra in \cite{J. L¨¹} with the coalgebra structures:
\begin{eqnarray*}
& \Delta(g)=g\otimes g,\ \Delta(x)=x\otimes 1+1\otimes x,\\
&\varepsilon(g)=1,\ \varepsilon(x)=0,
\end{eqnarray*}
and with the antipode
$$S(g)=g^{-1},\ S(x)=-x.$$

Moreover, we easily see that $(A,\Delta^\prime,\varepsilon,T)$ is a Hopf algebra with the following coalgebra structures:
\begin{eqnarray*}
& \Delta^\prime(g)=g\otimes g,\ \Delta^\prime(x)=x\otimes 1+g\otimes x,\\
&\varepsilon(g)=1,\ \varepsilon(x)=0,
\end{eqnarray*}
and with the antipode:
$$T(g)=g^{-1},\ T(x)=-g^{-1}x.$$

In the following, we prove taht $(A,\Delta,\Delta^\prime)$ is a Hopf brace.

In fact, we need to prove the compatibility condition (2.1) is satisfied. According to the facts $$
x_{1^\prime}\otimes x_{2^\prime}=x\otimes 1+g\otimes x,$$
$$
x_1\otimes x_2\otimes x_3=x\otimes 1\otimes 1+1\otimes x\otimes 1+1\otimes 1\otimes x,
$$ we have
\begin{eqnarray*}
x_{1^\prime}\otimes x_{2^\prime1}\otimes x_{2^\prime2}&=&x\otimes 1\otimes 1+g\otimes x_1\otimes x_2\\
&=&x\otimes 1\otimes 1+g\otimes x\otimes 1+g\otimes 1\otimes x,\\
x_{11^\prime}S(x_2)x_{31^\prime}\otimes x_{12^\prime}\otimes x_{32^\prime}&=&\underbrace{x_{1^\prime}S(1)1_{1^\prime}\otimes x_{2^\prime}\otimes 1_{2^\prime}}+\underbrace{1_{1^\prime}S(x)1_{1^\prime}\otimes 1_{2^\prime}\otimes 1_{2^\prime}}\\
&+&\underbrace{1_{1^\prime}S(1)x_{1^\prime}\otimes 1_{2^\prime}\otimes x_{2^\prime}}\\
&=&x\otimes 1\otimes 1+g\otimes x\otimes 1-x\otimes 1\otimes 1+x\otimes 1\otimes 1\\
&+&g\otimes 1\otimes x\\
&=&x\otimes 1\otimes 1+g\otimes x\otimes 1+g\otimes 1\otimes x,
 \end{eqnarray*}
so, we have
$$x_{1^\prime}\otimes x_{2^\prime1}\otimes x_{2^\prime2}=x_{11^\prime}S(x_2)x_{31^\prime}\otimes x_{12^\prime}\otimes x_{32^\prime}.
$$

In a similar way, we can prove that the other elements in the above equalities are satisfied. Hence $(A,\Delta,\Delta^\prime)$ is a Hopf brace.

(3) Let $H$ be a Hopf algebra. If $R=R^{\prime}_i\otimes R^{\prime\prime}_i\in H\o H$ is a normalized Harrison 2-cocycle in \cite{Nakajima1975}, that is, $R$ is satisfied the following conditions:
\begin{eqnarray*}
r^{\prime}_iR^{\prime}_{i1}\otimes r^{\prime\prime}_iR^{\prime}_{i2}\otimes R^{\prime\prime}_i&=&R^{\prime}_{i}\otimes
r^{\prime}_iR^{\prime\prime}_{i1}\otimes r^{\prime\prime}_iR^{\prime\prime}_{i2},\\
\varepsilon(R^{\prime}_i)R^{\prime\prime}_i&=&1\ =\ R^{\prime}_{i}\varepsilon(R^{\prime\prime}_i),
 \end{eqnarray*}
where $r$ is a copy of $R$. Define a comultiplication on $H$ as follows if $R$ is invertible with the inverse $R^{-1}=R^{\prime-1}_i\otimes R^{\prime\prime-1}_i:$
$$
\Delta_R(h)=R^{\prime}_ih_1R^{\prime-1}_i\otimes R^{\prime\prime}_ih_2R^{\prime\prime-1}_i
$$
for any $h\in H$. Then,  $(H,\D_R,S^R)$ is a Hopf algebra with the same counit, where $S^{R}(x)=R^{\prime}_iS(R^{\prime\prime}_i)S(x)S(R^{\prime-1}_i)R^{\prime\prime-1}_i$.

Let $(H, R)$ be a Long copaired Hopf algebra in \cite{Zhang2006}, that is, there is an invertible element $R=R^{\prime}_i\otimes R^{\prime\prime}_i\in H\o H$ such that the following conditions are satisfied:

$(LC1)\ \ R^{\prime}_ix\o R^{\prime\prime}_i=xR^{\prime}_i\o R^{\prime\prime}_i,$ for any $x\in H,$

$(LC2)\ \ \varepsilon(R^{\prime}_i)R^{\prime\prime}_i=1,$

$(LC3)\ \ R^{\prime}_{i1}\o R^{\prime}_{i2}\o R^{\prime\prime}_i=R^{\prime}_i\o r^{\prime}_i\o r^{\prime\prime}_iR^{\prime\prime}_i,$

$(LC4)\ \ R^{\prime}_i\varepsilon(R^{\prime\prime}_i)=1,$

$(LC5)\ \ R^{\prime}_i\o R^{\prime\prime}_{i1}\o R^{\prime\prime}_{i2}=R^{\prime}_ir^{\prime}_i\o R^{\prime\prime}_i\o r^{\prime\prime}_i.$

Suppose that $(H, R)$ is a Long copaired Hopf algebra. Then, according to Example 3.1 in \cite{Y. Chen}, we know that $R$ is a normalized Harrison 2-cocycle, so, we obtain a new Hopf algebra $(H, \Delta_R, S^R)$, and hence $(H, \Delta, \Delta_R)$ is a Hopf brace.

Indeed, we only need to check that the condition (2.1) is satisfied: for any $h\in H$,
\begin{eqnarray*}
LHB&=&R^{\prime}_ih_1R^{\prime-1}_i\o (R^{\prime\prime}_ih_2R^{\prime\prime-1}_i)_1\o (R^{\prime\prime}_ih_2R^{\prime\prime-1}_i)_2\\
&=&R^{\prime}_ih_1R^{\prime-1}_i\o R^{\prime\prime}_{i1}h_2R^{\prime\prime-1}_{i1}\o R^{\prime\prime}_{i2}h_3R^{\prime\prime-1}_{i2}\\
&\stackrel{(LC5)}=&R^{\prime}_ir^{\prime}_ih_1R^{\prime-1}_i\o R^{\prime\prime}_ih_2R^{\prime\prime-1}_{i1}\o r^{\prime\prime}_ih_3R^{\prime\prime-1}_{i2}\\
&=&R^{\prime}_ir^{\prime}_ih_1r^{\prime-1}_iR^{\prime-1}_i\o R^{\prime\prime}_ih_2R^{\prime\prime-1}_i\o r^{\prime\prime}_ih_3r^{\prime\prime-1}_i,
\end{eqnarray*}
where $R^{\prime-1}_i\o R^{\prime\prime-1}_{i1}\o R^{\prime\prime-1}_{i2}=r^{\prime-1}_iR^{\prime-1}_i\o R^{\prime\prime-1}_i\o R^{\prime\prime-1}_i$ by $(LC5)$.
\begin{eqnarray*}
RHB&=&R^{\prime}_ih_1R^{\prime-1}_iS(h_3)r^{\prime}_ih_4r^{\prime-1}_i\o R^{\prime\prime}_ih_2R^{\prime\prime-1}_i\o r^{\prime\prime}_ih_5r^{\prime\prime-1}_i\\
&\stackrel{(LC1)}=&R^{\prime}_ih_1R^{\prime-1}_iS(h_3)h_4r^{\prime}_ir^{\prime-1}_i\o R^{\prime\prime}_ih_2R^{\prime\prime-1}_i\o r^{\prime\prime}_ih_5r^{\prime\prime-1}_i\\
&=&R^{\prime}_ih_1R^{\prime-1}_ir^{\prime}_ir^{\prime-1}_i\o R^{\prime\prime}_ih_2R^{\prime\prime-1}_i\o r^{\prime\prime}_ih_3r^{\prime\prime-1}_i\\
&=&R^{\prime}_ir^{\prime}_ih_1R^{\prime-1}_ir^{\prime-1}_i\o R^{\prime\prime}_ih_2R^{\prime\prime-1}_i\o r^{\prime\prime}_ih_3r^{\prime\prime-1}_i\\
&=&R^{\prime}_ir^{\prime}_ir^{\prime-1}_ih_1R^{\prime-1}_i\o R^{\prime\prime}_ih_2R^{\prime\prime-1}_i\o r^{\prime\prime}_ih_3r^{\prime\prime-1}_i\\
&=&R^{\prime}_ir^{\prime}_ih_1r^{\prime-1}_iR^{\prime-1}_i\o R^{\prime\prime}_ih_2R^{\prime\prime-1}_i\o r^{\prime\prime}_ih_3r^{\prime\prime-1}_i,\\
\end{eqnarray*}
where $R^{\prime-1}_ix\o R^{\prime\prime-1}_i=xR^{\prime-1}_i\o R^{\prime\prime-1}_i$ by $(LC1)$, for any $x\in H$. So, (2.1) is satisfied, and hence $(H, \Delta, \Delta_R)$ is a Hopf brace.

\vspace{2mm}

Let $(H,\D,\D^\prime)$ and $(G,\D,\D^\prime)$ be Hopf braces. A homomorphism of Hopf braces $f: (H,\D,\D^\prime)\rightarrow (G,\D,\D^\prime)$ is a linear map $f$ such that $f: H_{\D}\rightarrow G_{\D}$ and $f: H_{\D^\prime}\rightarrow G_{\D^\prime}$ are Hopf algebra homomorphism. It is easy to see that Hopf braces form a category.

Fix a Hopf algebra $(H, m, 1, \Delta, \varepsilon, S)$. Let $\mathcal{CB}(H)$ be the full subcategory of the category of Hopf braces with objects $(H, \D, \D^\prime)$. This means that the objects of $\mathcal{CB}(H)$ are Hopf braces such that the first Hopf algebra structure is that of $H_{\D}$.

\vspace{2mm}

\noindent {\bf Lemma 2.4}\ \ Let $(H,\Delta,\Delta^\prime)$ be a Hopf brace. Then
$$
S(h_1)_{1^\prime}h_2\otimes S(h_1)_{2^\prime}=S(h_1)h_{21^\prime}\otimes S(h_{22^\prime})\eqno(2.2)
$$
for any $h\in H.$

\noindent {\bf Proof.}\ \ According to Eq.(2.1), we obtain the following equation:
$$
h\otimes 1=h_{1^\prime}\otimes h_{2^\prime 1}S(h_{2^\prime 2})=h_{11^\prime}S(h_2)h_{31^\prime}\otimes h_{12^\prime}S(h_{32^\prime})\eqno(2.3)
$$
for any $h\in H.$

So, we have
\begin{eqnarray*}
S(h_1)_{1^\prime}h_2\otimes S(h_1)_{2^\prime}&\stackrel{(2.3)}{=}&S(h_1)_{1^\prime}h_{21^\prime}S(h_3)h_{41^\prime}\otimes S(h_1)_{2^\prime}h_{22^\prime}S(h_{42^\prime})\\
&=&(S(h_1)h_2)_{1^\prime}S(h_3)h_{41^\prime}\otimes (S(h_1)h_2)_{2^\prime}S(h_{42^\prime})\\
&=&S(h_1)h_{21^\prime}\otimes S(h_{22^\prime}).
 \end{eqnarray*}
\hfill $\square$
\vspace{2mm}

\noindent {\bf Lemma 2.5}\ \ Let $(H,\Delta,\Delta^\prime)$ is a Hopf brace. Then, $H$ is a left $(H,\Delta^\prime)$-comodule coalgebra with
$$
\rho(h)=S(h_1)h_{21^\prime}\otimes h_{22^\prime}
$$
for any $h\in H.$

\noindent {\bf Proof.}\ \ Firstly, we can check that $H$ is a left $(H,\Delta^\prime)$-comodule: for any $h\in H,$ we have
\begin{eqnarray*}
(id\otimes\rho)\rho(h)&=&S(h_1)h_{21^\prime}\otimes \rho(h_{22^\prime})\\
&=&S(h_1)h_{21^\prime}\otimes S(h_{22^\prime 1})h_{22^\prime 21^\prime}\otimes h_{22^\prime 22^\prime}\\
&\stackrel{(2.1)}{=}&S(h_1)h_{21^\prime} S(h_3)h_{41^\prime}\otimes S(h_{22^\prime}) h_{42^\prime} \otimes h_{43^\prime}\\
&\stackrel{(2.2)}{=}&S(h_1)_{1^\prime}h_2S(h_3)h_{41^\prime}\otimes S(h_1)_{2^\prime}h_{42^\prime}\otimes h_{43^\prime}\\
&=&S(h_1)_{1^\prime}h_{21^\prime}\otimes S(h_1)_{2^\prime}h_{22^\prime}\otimes h_{23^\prime}\\
&=&(S(h_1)h_{21^\prime})_{1^\prime}\otimes (S(h_1)h_{21^\prime})_{2^\prime}\otimes h_{22^\prime}\\
&=&(\Delta^\prime\otimes id)(S(h_1)h_{21^\prime}\otimes h_{22^\prime})\\
&=&(\Delta^\prime\otimes id)\rho(h),\\
(\varepsilon\otimes id)\rho(h)&=&h.
 \end{eqnarray*}

So, it is a left $(H,\Delta^\prime)$-comodule. In addition, we know
\begin{eqnarray*}
h_{(-1)}\otimes\Delta(h_{(0)})&=&S(h_1)h_{21^\prime}\otimes h_{22^\prime 1}\otimes h_{22^\prime2}\\
&\stackrel{(2.1)}{=}&S(h_1)h_{211^\prime}S(h_{22})h_{231^\prime}\otimes h_{212^\prime}\otimes h_{232^\prime}\\
&=&S(h_1)h_{21^\prime}S(h_3)h_{41^\prime}\otimes h_{22^\prime}\otimes h_{42^\prime}\\
&=&h_{1(-1)}h_{2(-1)}\otimes h_{1(0)}\otimes h_{2(0)},\\
h_{(-1)}\varepsilon(h_{(0)})&=&S(h_1)h_{21^\prime}\varepsilon(h_{22^\prime})\\
&=&S(h_1)h_2=\varepsilon(h)1_H.
 \end{eqnarray*}

The proof is complete.
\hfill $\square$

\vspace{2mm}

\noindent {\bf Remark 2.6}\ \ It follows from the Lemma 2.4 that
$$h_{1^\prime}\otimes h_{2^\prime}=h_1h_{2(-1)}\o h_{2(0)},\eqno{(2.4)}$$
$$h_1\otimes h_2=h_{1^\prime}T(h_{2^\prime(-1)})\o h_{2^\prime(0)},\eqno{(2.5)}
$$
for any $h\in H.$

\noindent {\bf Proof.} In fact, for any $h\in H$, we have
 \begin{eqnarray*}
h_1h_{2(-1)}\o h_{2(0)}&=&h_1S(h_2)h_{31^\prime}\otimes h_{32^\prime}\\
&=&\varepsilon(h_1)h_{21^\prime}\otimes h_{22^\prime}=h_{1^\prime}\otimes h_{2^\prime},\\
h_{1^\prime}T(h_{2^\prime(-1)})\o h_{2^\prime(0)}&=&h_{1^\prime}T(S(h_{2^\prime1})h_{2^\prime21^\prime})\otimes h_{2^\prime22^\prime}\\
&\stackrel{(2.1)}{=}&h_{11^\prime}S(h_2)h_{31^\prime}T(S(h_{12^\prime})h_{32^\prime})\otimes h_{33^\prime}\\
&=&h_{11^\prime}S(h_2)h_{31^\prime}T(h_{32^\prime})T(S(h_{12^\prime}))\otimes h_{33^\prime}\\
&=&h_{11^\prime}S(h_2)\varepsilon(h_{31^\prime})T(S(h_{12^\prime}))\otimes h_{32^\prime}\\
&=&h_{11^\prime}S(h_2)T(S(h_{12^\prime}))\otimes h_3\\
&\stackrel{(2.4)}{=}&h_1S(h_2)h_{31^\prime}S(h_4)T(S(h_{32^\prime}))\otimes h_5\\
&\stackrel{(2.2)}{=}&h_1S(h_2)_{1^\prime}h_3S(h_4)T(S(h_2)_{2^\prime})\otimes h_5\\
&=&h_1S(h_2)_{1^\prime}T(S(h_2)_{2^\prime})\otimes h_3\\
&=&h_1\otimes h_2.
\end{eqnarray*}
\hfill $\square$

\vspace{2mm}

\noindent {\bf Definition 2.7}\ \ A Hopf brace $(A,\Delta,\Delta^\prime)$ is said to be commutative if the underlying algebra $A$ is commutative.

\vspace{2mm}

\noindent{\bf Lemma 2.8}\ \ Let $(A,\Delta,\Delta^\prime)$ be a commutative Hopf brace. Then the following conclusions hold:

(1) $A$ is a left $A_{\Delta^\prime}$-comodule algebra via $$\rho(a)\equiv a_{(-1)}\otimes a_{(0)}=S(a_1)a_{21^\prime}\otimes a_{22^\prime},$$ for any $a\in A$.

(2) $A$ is a right $A_{\Delta^\prime}$-comodule algebra via
$$\varphi(a)\equiv a_{[0]}\otimes a_{[1]}=T(a_{1^\prime})_{(-1)}a_{2^\prime}\o T(a_{1^\prime})_{(0)}a_{3^\prime}=S(T(a_{1^\prime})_1)T(a_{1^\prime})_{21^\prime}a_{2^\prime}\otimes T(a_{1^\prime})_{22^\prime}a_{3^\prime},$$
for any $a\in A$.

(3) $(id\otimes S)\rho(a)=\rho(S(a))$, for all $a \in A$. That is, $a_{(-1)}\otimes S(a_{(0)})=S(a)_{(-1)}\otimes S(a)_{(0)}$.

\noindent{\bf Proof.}
(1) By Lemma 2.5, we know $A$ is a left $A_{\Delta^\prime}$-comodule. In the following, we only need to prove that $\rho$ is an algebra map.

In fact, for any $a,b\in A$, we have
\begin{eqnarray*}
\rho(ab)&=&S((ab)_1)(ab)_{21^\prime}\otimes(ab)_{22^\prime}\\
&=&S(a_1)S(b_1)a_{21^\prime}b_{21^\prime}\otimes a_{22^\prime}b_{22^\prime}\\
&=&S(a_1)a_{21^\prime}S(b_1)b_{21^\prime}\otimes a_{22^\prime}b_{22^\prime}\\
&=&\rho(a)\rho(b).
\end{eqnarray*}

It is obvious that $\rho(1_A)=1_A\o 1_A$. So, $A$ is a left $A_{\Delta^\prime}$-comodule algebra.

(2) Firstly, $A$ is a right $A_{\Delta^\prime}$-comodule, because for any $a\in A$,
\begin{eqnarray*}
(\varphi\otimes id)\varphi(a)&=&a_{[0][0]}\otimes a_{[0][1]}\otimes a_{[1]}\\
&=&(T(a_{1^\prime})_{(-1)}a_{2^\prime})_{[0]}\otimes (T(a_{1^\prime})_{(-1)}a_{2^\prime})_{[1]}\otimes T(a_{1^\prime})_{(0)}a_{3^\prime}\\
&=&T(T(a_{1^\prime})_{(-1)1^\prime} a_{2^\prime})_{(-1)}T(a_{1^\prime})_{(-1)2^\prime}a_{3^\prime}\otimes \\ &&T(T(a_{1^\prime})_{(-1)1^\prime}a_{2^\prime})_{(0)}T(a_{1^\prime})_{(-1)3^\prime}a_{4^\prime}\otimes T(a_{1^\prime})_{(0)}a_{5^\prime}\\
&=&T(T(a_{1^\prime})_{(-1)1^\prime})_{(-1)}T(a_{1^\prime})_{(-1)2^\prime}T(a_{2^\prime})_{(-1)}a_{3^\prime}\\
&&\otimes T(T(a_{1^\prime})_{(-1)1^\prime})_{(0)}T(a_{1^\prime})_{(-1)3^\prime}T(a_{2^\prime})_{(0)}a_{4^\prime}\otimes T(a_{1^\prime})_{(0)}a_{5^\prime}\\
&=&T(T(a_{1^\prime})_{(-1)1^\prime})_{(-1)}T(a_{1^\prime})_{(-1)2^\prime1^\prime}T(a_{2^\prime})_{(-1)}a_{3^\prime}\\
&&\otimes T(T(a_{1^\prime})_{(-1)1^\prime})_{(0)}T(a_{1^\prime})_{(-1)2^\prime2^\prime}T(a_{2^\prime})_{(0)}a_{4^\prime}\otimes T(a_{1^\prime})_{(0)}a_{5^\prime}\\
&=&T(T(a_{1^\prime})_{(-1)})_{(-1)}T(a_{1^\prime})_{(0)(-1)1^\prime}T(a_{2^\prime})_{(-1)}a_{3^\prime}\\
&&\otimes T(T(a_{1^\prime})_{(-1)})_{(0)}T(a_{1^\prime})_{(0)(-1)2^\prime}T(a_{2^\prime})_{(0)}a_{4^\prime}\otimes T(a_{1^\prime})_{(0)(0)}a_{5^\prime}\\
&=&T(T(a_{1^\prime})_{2^\prime(-1)})_{(-1)}T(a_{1^\prime})_{2^\prime(0)(-1)1^\prime}T(a_{1^\prime})_{1^\prime(-1)}a_{2^\prime}\\
&&\otimes T(T(a_{1^\prime})_{2^\prime(-1)})_{(0)}T(a_{1^\prime})_{2^\prime(0)(-1)2^\prime}T(a_{1^\prime})_{1^\prime(0)}a_{3^\prime}\otimes T(a_{1^\prime})_{2^\prime(0)(0)}a_{4^\prime}\\
&=&(T(a_{1^\prime})_{1^\prime}T(T(a_{1^\prime})_{2^\prime(-1)}))_{(-1)}T(a_{1^\prime})_{2^\prime(0)(-1)1^\prime}a_{2^\prime}\\
&&\otimes (T(a_{1^\prime})_{1^\prime}T(T(a_{1^\prime})_{2^\prime(-1)}))_{(0)}T(a_{1^\prime})_{2^\prime(0)(-1)2^\prime}a_{3^\prime}\otimes T(a_{1^\prime})_{2^\prime(0)(0)}a_{4^\prime}\\
&\stackrel{(2.5)}{=}&T(a_{1^\prime})_{1(-1)}T(a_{1^\prime})_{2(-1)1^\prime}a_{2^\prime}\otimes T(a_{1^\prime})_{1(0)}T(a_{1^\prime})_{2(-1)2^\prime}a_{3^\prime}\otimes T(a_{1^\prime})_{2(0)}a_{4^\prime}\\
&=&S(T(a_{1^\prime})_1)T(a_{1^\prime})_{21^\prime}S(T(a_{1^\prime})_3)_{1^\prime}T(a_{1^\prime})_{41^\prime}a_{2^\prime}\\
&&\otimes T(a_{1^\prime})_{22^\prime}S(T(a_{1^\prime})_{3})_{2^\prime}T(a_{1^\prime})_{42^\prime}a_{3^\prime}\otimes T(a_{1^\prime})_{43^\prime}a_{4^\prime}\\
&=&S(T(a_{1^\prime})_1)T(a_{1^\prime})_{21^\prime}a_{2^\prime}\otimes T(a_{1^\prime})_{22^\prime}a_{3^\prime}\otimes T(a_{1^\prime})_{23^\prime}a_{4^\prime}\\
&=&T(a_{1^\prime})_{(-1)}a_{2^\prime}\otimes T(a_{1^\prime})_{(0)1^\prime}a_{3^\prime}\otimes T(a_{1^\prime})_{(0)2^\prime}a_{4^\prime}\\
&=&a_{[0]}\otimes a_{[1]1^\prime}\otimes a_{[1]2^\prime}\\
&=&(id\otimes \Delta^\prime)\varphi(a),\\
\varphi(a)\varphi(b)&=&(a_{[0]}\otimes a_{[1]})(b_{[0]}\otimes b_{[1]}).
\end{eqnarray*}

Secondly, since $(A,\Delta,\Delta^\prime)$ is commutative , we easily see that $\varphi$ is an algebra map. Thus, $A$ is a right $A_{\Delta^\prime}$-comodule algebra.

(3) For any $a\in A$, by Lemma 2.5, we have
\begin{eqnarray*}
S(a)_{(-1)}\otimes S(a)_{(0)}&=a_2S(a_1)_{1^\prime}\otimes S(a_1)_{2^\prime}=S(a_1)_{1^\prime}a_2\otimes S(a_1)_{2^\prime}\\
&=S(a_1)a_{21^\prime}\otimes S(a_{22^\prime})=a_{(-1)}\otimes S(a_{(0)}).
\end{eqnarray*} So, (3) is proved.
\hfill $\square$

\vspace{2mm}

\noindent{\bf Proposition 2.9}\ \ Let $(A,\Delta,\Delta^\prime)$ be a commutative Hopf brace. Then, the given map
$$c: A\o A\rightarrow A\o A, \ \
c(x\o y)=x_{(-1)}y_{[0]}\o x_{(0)}y_{[1]},
$$
 is a solution of the braid equation.

\noindent{\bf Proof.}\ \ Let $\gamma: A\o A\rightarrow A\o A$ be given by $\g(x\o y)=xy_{(-1)}\o y_{(0)}$. Then, it is easy to prove that $\g$ is invertible with the inverse $\g^{-1}(x\o y)=xT(y_{(-1)})\o y_{(0)}$.

Let $\s: A\o A\rightarrow A\o A$, $\s(x\o y)=y_2\o xS(y_1)y_3$.
By using Lemma 2.8 and Remark 2.6, we have
\begin{eqnarray*}
\g^{-1}c\g(x\o y)&=&\g^{-1}c(xy_{(-1)}\o y_{(0)})\\
&=&\g^{-1}((xy_{(-1)})_{(-1)}y_{(0)[0]}\o (xy_{(-1)})_{(0)}y_{(0)[1]})\\
&=&\g^{-1}(x_{(-1)}y_{(-1)(-1)}T(y_{(0)1^\prime})_{(-1)}y_{(0)2^\prime}\o x_{(0)}y_{(-1)(0)}T(y_{(0)1^\prime})_{(0)}y_{(0)3^\prime})\\
&=&\g^{-1}(S(x_1)x_{21^\prime}(S(y_1)y_{21^\prime})_{(-1)}T(y_{22^\prime1^\prime})_{(-1)}y_{22^\prime2^\prime}\\
&&\o x_{22^\prime} (S(y_1)y_{21^\prime})_{(0)}T(y_{22^\prime1^\prime})_{(0)}y_{22^\prime3^\prime})\\
&=&\g^{-1}(S(x_1)x_{21^\prime}S(y_1)_{(-1)}y_{21^\prime(-1)}T(y_{22^\prime1^\prime})_{(-1)}y_{22^\prime2^\prime}\\
&&\o x_{22^\prime} S(y_1)_{(0)}y_{21^\prime(0)}T(y_{22^\prime1^\prime})_{(0)}y_{22^\prime3^\prime})\\
&=&\g^{-1}(S(x_1)x_{21^\prime}y_{1(-1)}y_{21^\prime(-1)}T(y_{22^\prime1^\prime})_{(-1)}y_{22^\prime2^\prime}\\
&&\o x_{22^\prime} S(y_{1(0)})y_{21^\prime(0)}T(y_{22^\prime1^\prime})_{(0)}y_{22^\prime3^\prime})\\
&=&S(x_1)x_{21^\prime}y_{1(-1)}y_{21^\prime(-1)}T(y_{22^\prime1^\prime})_{(-1)}y_{22^\prime2^\prime}T(x_{22^\prime(-1)})T(S(y_{1(0)})_{(-1)})\\
&&T(y_{21^\prime(0)(-1)})T(T(y_{22^\prime1^\prime})_{(0)(-1)})T(y_{22^\prime3^\prime(-1)}) \\
&&\o x_{22^\prime(0)} S(y_{1(0)})_{(0)}y_{21^\prime(0)(0)}T(y_{22^\prime1^\prime})_{(0)(0)}y_{22^\prime3^\prime(0)}\\
&=&S(x_1)x_{21^\prime}T(x_{22^\prime(-1)})y_{1(-1)}y_{21^\prime(-1)}T(y_{22^\prime1^\prime})_{(-1)}y_{22^\prime2^\prime}T(S(y_{1(0)})_{(-1)})\\
&&T(y_{21^\prime(0)(-1)})T(T(y_{22^\prime1^\prime})_{(0)(-1)})T(y_{22^\prime3^\prime(-1)}) \\
&&\o x_{22^\prime(0)} S(y_{1(0)})_{(0)}y_{21^\prime(0)(0)}T(y_{22^\prime1^\prime})_{(0)(0)}y_{22^\prime3^\prime(0)}\\
&=&S(x_1)x_2y_{1(-1)}y_{21^\prime(-1)}T(y_{22^\prime1^\prime})_{(-1)}y_{22^\prime2^\prime}T(S(y_{1(0)})_{(-1)})T(y_{21^\prime(0)(-1)})\\
&&T(T(y_{22^\prime1^\prime})_{(0)(-1)})T(y_{22^\prime3^\prime(-1)})\o x_3 S(y_{1(0)})_{(0)}y_{21^\prime(0)(0)}T(y_{22^\prime1^\prime})_{(0)(0)}y_{22^\prime3^\prime(0)}\\
&=&y_{1(-1)}y_{21^\prime(-1)}T(y_{22^\prime1^\prime})_{(-1)}y_{22^\prime2^\prime}T(y_{1(0)(-1)})T(y_{21^\prime(0)(-1)})
T(T(y_{22^\prime1^\prime})_{(0)(-1)})\\
&&T(y_{22^\prime3^\prime(-1)}) \o x S(y_{1(0)(0)})y_{21^\prime(0)(0)}T(y_{22^\prime1^\prime})_{(0)(0)}y_{22^\prime3^\prime(0)}\\
&=&y_{1(-1)1^\prime}y_{21^\prime(-1)1^\prime}T(y_{22^\prime})_{(-1)1^\prime}y_{23^\prime}T(y_{1(-1)2^\prime})T(y_{21^\prime(-1)2^\prime})
T(T(y_{22^\prime})_{(-1)2^\prime})\\
&&T(y_{24^\prime(-1)}) \o x S(y_{1(0)})y_{21^\prime(0)}T(y_{22^\prime})_{(0)}y_{24^\prime(0)}\\
\end{eqnarray*}
\begin{eqnarray*}
&=&y_{1(-1)1^\prime}T(y_{1(-1)2^\prime})y_{21^\prime(-1)1^\prime}T(y_{21^\prime(-1)2^\prime})T(y_{22^\prime})_{(-1)1^\prime}T(T(y_{22^\prime})_{(-1)2^\prime})
y_{23^\prime}\\
&&T(y_{24^\prime(-1)}) \o x S(y_{1(0)})y_{21^\prime(0)}T(y_{22^\prime})_{(0)}y_{24^\prime(0)}\\
&=&\v(y_{1(-1)})\v(y_{21^\prime(-1)})\v(T(y_{22^\prime})_{(-1)})y_{23^\prime}T(y_{24^\prime(-1)})\\
 &&\o x S(y_{1(0)})y_{21^\prime(0)}T(y_{22^\prime})_{(0)}y_{24^\prime(0)}\\
&=&y_{23^\prime}T(y_{24^\prime(-1)}) \o x S(y_{1})y_{21^\prime}T(y_{22^\prime})y_{24^\prime(0)}\\
&=&y_{22^\prime}T(y_{23^\prime(-1)}) \o x S(y_{1})\v(y_{21^\prime})y_{23^\prime(0)}\\
&=&y_{21^\prime}T(y_{22^\prime(-1)}) \o x S(y_{1})y_{22^\prime(0)}=y_{21} \o x S(y_{1})y_{22}\\
&=&y_{2} \o x S(y_{1})y_{3}\\
&=&\s(x\o y),
\end{eqnarray*}
for any $x,y\in A$. From this it follows that $c$ is an invertible solution of the braid equation.
\hfill $\square$

\vspace{2mm}

\noindent {\bf Definition 2.10}\ \
Let $H$ and $A$ be Hopf algebras. Assume that $A$ be an $H$-comodule coalgebra. A {\it bijective 1-cocycle} is an algebra isomorphism $\pi:A \rightarrow H$ such that
$$
\pi(a)_1\otimes\pi(a)_2=\pi(a_1)a_{2(-1)}\otimes \pi(a_{2(0)}),\eqno(2.6)
$$ for any $a \in A.$

\vspace{2mm}

\noindent {\bf Remark 2.11}\ \
(1) Any bijective 1-cocycle $\pi$ satisfies $\v_H\pi=\v_A$. Indeed, applying $\v_H\o id$ to Eq.(2.6), we have $\pi(a)=\v_H(\pi(a_1))\pi(a_2)$. Since $\pi$ is a bijection, we get $a=\v_H(\pi(a_1))a_2$. Applying $\v_A$ to this equation, we obtain $\v_H\pi=\v_A$.

(2) Let $\pi:A \xrightarrow[]{} H$ and $\eta:B \xrightarrow[]{} K$ be two bijective 1-cocycles. A morphism between these bijective 1-cocycles is a pair $(f, g)$ of Hopf algebra maps $f:K \xrightarrow[]{} H,~g:B \xrightarrow[]{} A,$ such that the following conditions are satisfied:
$$\begin{array}{rllr}
&\pi g=f\eta,\\
&g(b)_{(-1)}\otimes g(b)_{(0)}=f(b_{(-1)})\otimes g(b_{(0)}),
\end{array}$$ for any $b \in B.$ It is easy to see that bijective 1-cocycles form a category. Fix a Hopf algebra $A$, we assume that $\mathcal{C}(A)$ is the full subcategory of the category of bijective 1-cocycles with objects $\pi:A \xrightarrow[]{} H.$

\vspace{2mm}

\noindent {\bf Theorem 2.12}\ \
Let $A$ be a Hopf algebra. Then, the category $\mathcal{CB}(A)$ of Hopf braces is equivalent to the full subcategory $\mathcal{C}(A)$ of of bijective 1-cocycles.\\
\noindent{\bf Proof.}\ \
We claim that $F:\mathcal{CB}(A) \xrightarrow[]{} \mathcal{C}(A)$ is a functor given by
$$
F(A,\Delta,\Delta^\prime)=(id_A:A_{\Delta} \xrightarrow[]{} A_{\Delta^\prime}),~F(f)=(f,f)
$$ for the morphism $f:(A,\Delta,\Delta^\prime) \xrightarrow[]{} (A,\Delta,\Delta^\circledast)$.

We prove that $\pi=id_A:A_{\Delta} \rightarrow A_{\Delta^\prime}$ is a bijective 1-cocycle. By Lemma 2.5, $A_{\Delta}$ is a $A_{\Delta^\prime}$-comodule coalgebra and for any $a \in A,$
\begin{eqnarray*}
\pi(a_1)a_{2(-1)}\otimes\pi(a_{2(0)})&=&\pi(a_1)S(a_2)a_{31^\prime}\otimes\pi(a_{32^\prime})\\
&=&a_1S(a_2)a_{31^\prime}\otimes a_{32^\prime}\\
&=&\Delta^\prime(a)=\pi(a)_{1^\prime}\otimes \pi(a)_{2^\prime}.
\end{eqnarray*}

Now $(f,f)$ is a morphism of bijective 1-cocycles since $f$ is a Hopf algebra map from $(A,\Delta)$ to $(A,\Delta^\prime)$. Hence the claim follows.

Now we define a functor $G:\mathcal{C}(A)\rightarrow \mathcal{CB}(A)$ as follows. First, $G(\pi:A \rightarrow H)=(A,\Delta,\Delta^\prime)$, where the new comultiplication is given by
$$\begin{array}{rllr}
\Delta^\prime(a)=\pi^{-1}(\pi(a)_1)\otimes \pi^{-1}(\pi(a)_2)
\end{array}$$ for any $a \in A.$

In what follows, we prove that $(A,\Delta,\Delta^\prime)$ is a Hopf brace. It is easy to check that $(A, m, 1, \Delta^\prime,\varepsilon, T\equiv \pi^{-1}S\pi)$ is a Hopf algebra.

To prove that $(A,\Delta,\Delta^\prime)$ is a Hopf brace, for any $a \in A$, we have
\begin{eqnarray*}
&&a_{11^\prime}S(a_2)a_{31^\prime}\otimes a_{12^\prime}\otimes a_{32^\prime}\\
&=&\pi^{-1}(\pi(a_1)_1)S(a_2)\pi^{-1}(\pi(a_3)_1)\otimes\pi^{-1}(\pi(a_1)_2)\otimes\pi^{-1}(\pi(a_3)_2)\\
&\stackrel{(2.6)}{=}&\pi^{-1}(\pi(a_{11})a_{12(-1)})S(a_2)\pi^{-1}(\pi(a_{31})a_{32(-1)})\otimes \pi^{-1}(\pi(a_{12(0)}))\otimes \pi^{-1}(\pi(a_{32(0)}))\\
&=&a_1\pi^{-1}(a_{2(-1)})S(a_3)a_4\pi^{-1}(a_{5(-1)})\otimes a_{2(0)}\otimes a_{5(0)}\\
&=&a_1\pi^{-1}(a_{2(-1)})\pi^{-1}(a_{3(-1)})\otimes a_{2(0)}\otimes a_{3(0)}\\
&=&a_1\pi^{-1}(a_{21(-1)})\pi^{-1}(a_{22(-1)})\otimes a_{21(0)}\otimes a_{22(0)}\\
&=&a_1\pi^{-1}(a_{2(-1)})\otimes a_{2(0)1}\otimes a_{2(0)2}\\
&=&\pi^{-1}(\pi(a_1)a_{2(-1)})\otimes \pi^{-1}(\pi(a_{2(0)}))_1\otimes \pi^{-1}(\pi(a_{2(0)}))_2\\
&\stackrel{(2.6)}{=}&\pi^{-1}(\pi(a)_1)\otimes\pi^{-1}(\pi(a)_2)_1\otimes\pi^{-1}(\pi(a)_2)_2\\
&=&a_{1^\prime}\otimes a_{2^\prime1}\otimes a_{2^\prime2}.
\end{eqnarray*}

Thus, $(A,\Delta,\Delta^\prime)$ is a Hopf brace.

For any given morphism $(f,g)$ between bijective 1-cocycles $\pi$ and $\eta$, we define $G(f,g)=g.$
Then, for any $b \in B$, we have
\begin{eqnarray*}
\Delta^\prime g(b)&=&\pi^{-1}(\pi(g(b))_1)\otimes\pi^{-1}(\pi(g(b))_2)\\
&=&\pi^{-1}(f(\eta(b))_1)\otimes\pi^{-1}(f(\eta(b))_2)\\
&=&\pi^{-1}(f(\eta(b)_1))\otimes\pi^{-1}(f(\eta(b)_2))\\
&=&g\eta^{-1}(\eta(b)_1)\otimes g\eta^{-1}(\eta(b)_2)\\
&=&(g\o g)(\eta^{-1}(\eta(b)_1)\otimes \eta^{-1}(\eta(b)_2))\\
&=&(g\o g)\Delta^\circledast(b).
\end{eqnarray*}

So, $G$ is a function. Clearly $G\circ F\simeq id_{\mathcal{CB}(A)}$ and $F\circ G\simeq id_{\mathcal{C}(A)}.$
\hfill $\square$

\vspace{3mm}

\begin{center}
{\bf \S3\quad Hopf brace and Hopf matched pair}
\end{center}

\vspace{3mm}

In this section, we mainly study structures of commutative Hopf braces, build a correspondence between Hopf braces and Hopf matched pairs, and prove that the category $\mathcal{CB}(A)$ of Hopf braces is equivalent to the category $\mathcal{M}(A)$ of Hopf matched pairs.

\vspace{2mm}

\textbf{Definition 3.1~}~ Let $A$ and $H$ be Hopf algebras. A {\it Hopf matched pair} is a pair $(A,H)$ with two coactions
\begin{eqnarray*}
H\stackrel{\varphi}\longrightarrow H\otimes A\stackrel{\rho}\longleftarrow A
\end{eqnarray*}
such that $(A, \rho)$ is a left $H$-comodule algebra, $(H,\varphi)$ a right $A$-comodule algebra, and the following compatibilities hold:

$(HM1)\ \ a_{(-1)}\v_A(a_{(0)})=\v_A(a)1_H,~~\v_H(h_{[0]})h_{[1]}=\v_H(h)1_A,$

$(HM2)\ \ a_{(-1)}\otimes a_{(0)1}\otimes a_{(0)2}=a_{1(-1)}a_{2(-1)[0]}\otimes a_{1(0)}a_{2(-1)[1]}\otimes a_{2(0)},$

$(HM3)\ \ h_{[0]1}\otimes h_{[0]2}\otimes h_{[1]}=h_{1[0]}\otimes h_{1[1](-1)}h_{2[0]}\otimes h_{1[1](0)}h_{2[1]},$

$(HM4)\ \ h_{[0]}a_{(-1)}\o h_{[1]}a_{(0)}=a_{(-1)}h_{[0]}\o a_{(0)}h_{[1]},$
\\
for any $a\in A, h\in H$, where $\rho(a)$ is denoted by $a_{(-1)}\o a_{(0)}$ and $\varphi(h)$ denoted by $h_{[0]}\o h_{[1]}$.

\vspace{2mm}

\textbf{Example 3.2~} (1) Let $A=k[g,g^{-1},x]$ be a Hopf algebra as in Example 2.3, and let $H=k[X,a^{\pm},b^{\pm}]$ a Hopf algebra in \cite{J. L¨¹} with the following structures:
\begin{eqnarray*}
\Delta(a)&=&a\otimes a,\ \Delta(b)=b\otimes b,\ \Delta(X)=X\otimes ab+ab\otimes X\\
\varepsilon(a)&=&\varepsilon(b)=1,\ \varepsilon(X)=0,\\
S(a)&=&a^{-1},\ S(b)=b^{-1},\ S(X)=-a^{-2}b^{-2}X.
\end{eqnarray*}

Then, it is easy to get a Hopf matched pair $(A,H,\rho,\varphi)$ with coactions as follows:
\begin{eqnarray*}
\rho(g)&=&1\otimes g,\ \rho(x)=a\otimes x,\\
\varphi(X)&=&X\otimes g,\ \varphi(a)=a\otimes 1,\ \varphi(b)=b\otimes 1.
\end{eqnarray*}

(2) Let $H$ and $A$ be Hopf algebras. An invertible element $R=R_i^{\prime}\otimes R_i^{\prime\prime}$ in $H\otimes A$ is called a weak $R$-matrix of $H$ and $A$ in \cite{H. X. Chen} if the following conditions are satisfied:

$(WM1)\ (\D\o id)(R)=R_i^{\prime}\o r_i^{\prime}\o R_i^{\prime\prime}r_i^{\prime\prime},$

$(WM2)\ (id\o \D)(R)=R_i^{\prime}r_i^{\prime}\o r_i^{\prime\prime}\o R_i^{\prime\prime},$\\
where $r=r_i^{\prime}\otimes r_i^{\prime\prime}$ is a copy of $R$.

Then, by Lemma 1.3 in \cite{H. X. Chen}, $(H, A, \rho, \varphi)$ is a Hopf matched pair with the coactions as follows:
\begin{eqnarray*}
&\rho: H\rightarrow A\o H, \rho(h)=\tau(R)(1\o h)\tau(R^{-1}),\\
&\varphi: A\rightarrow A\o H, \varphi(a)=\tau(R)(a\o 1)\tau(R^{-1}),
\end{eqnarray*}
where $\tau$ is the twisted map and $R^{-1}$ the inverse of $R$.

\vspace{2mm}

In the following, we will show that there is a correspondence between Hopf braces and Hopf matched pairs.

\vspace{2mm}

\noindent{\bf Proposition 3.3}\ \ Let $(A,\Delta,\Delta^\prime)$ be a commutative Hopf brace. Then, $(A_{\Delta^\prime},A_{\Delta^\prime})$ is a Hopf matched pair with coactions as follows:
\begin{eqnarray*}
&&\rho(a)=a_{(-1)}\otimes a_{(0)}=S(a_1)a_{21^\prime}\otimes a_{22^\prime},\\
&&\varphi(a)=a_{[0]}\otimes a_{[1]}=T(a_{1^\prime})_{(-1)}a_{2^\prime}\o T(a_{1^\prime})_{(0)}a_{3^\prime}= S(T(a_{1^\prime})_1)T(a_{1^\prime})_{21^\prime}a_{2^\prime}\otimes T(a_{1^\prime})_{22^\prime}a_{3^\prime},
\end{eqnarray*}
for any $a\in A$.

\noindent{\bf Proof.}\ \
By Lemma 2.8, we need to prove that Eqs.$(HM1)-(HM3)$ hold.

In fact, by the above coactions, we have
\begin{eqnarray*}
 a_{1^\prime(-1)}a_{2^\prime[0]}\otimes a_{1^\prime(0)}a_{2^\prime[1]}
 &=&a_{1^\prime(-1)}T(a_{2^\prime})_{(-1)}a_{3^\prime}\otimes a_{1^\prime(0)}T(a_{2^\prime})_{(0)}a_{4^\prime}\\
 &=&\varepsilon(a_{1^\prime})a_{2^\prime}\otimes a_{3^\prime}\\
 &=&a_{1^\prime}\otimes a_{2^\prime}
\end{eqnarray*}
for any $a\in A$. Again by Remark 2.6, we get
\begin{eqnarray*}
&&a_{1^\prime(-1)}a_{2^\prime(-1)[0]}\otimes a_{1^\prime(0)}a_{2^\prime(-1)[1]}\otimes a_{2^\prime(0)}\\
&=&a_{1(-1)}a_{2(-1)(-1)}a_{2(0)(-1)[0]}\otimes a_{1(0)}a_{2(-1)(0)}a_{2(0)(-1)[1]}\otimes a_{2(0)(0)}\\
&=&a_{1(-1)}a_{2(-1)1^\prime(-1)}a_{2(-1)2^\prime[0]}\otimes a_{1(0)}a_{2(-1)1^\prime(0)}a_{2(-1)2^\prime[1]}\otimes a_{2(0)}\\
&=&a_{1(-1)}a_{2(-1)1^\prime}\otimes a_{1(0)}a_{2(-1)2^\prime}\otimes a_{2(0)}\\
&=&a_{1(-1)}a_{2(-1)}\otimes a_{1(0)}a_{2(0)(-1)}\otimes a_{2(0)(0)}\\
&=&a_{(-1)}\otimes a_{(0)1}a_{(0)2(-1)}\otimes a_{(0)2(0)}\\
&=&a_{(-1)}\otimes a_{(0)1^\prime}\otimes a_{(0)2^\prime},
\end{eqnarray*}
and so Eq.$(HM2)$ holds. Moreover, we have
\begin{eqnarray*}
\ \ \ \ \ \ \  \ \ \ \  \ \ \ \  \ \ \ \ \ &&a_{1^\prime[0]}\otimes a_{1^\prime[1](-1)}a_{2^\prime[0]}\otimes a_{1^\prime[1](0)}a_{2^\prime[1]}\\
&=&T(a_{1^\prime})_{(-1)}a_{2^\prime}\otimes T(a_{1^\prime})_{(0)(-1)}a_{3^\prime(-1)}T(a_{4^\prime})_{(-1)}a_{5^\prime}\otimes T(a_{1^\prime})_{(0)(0)}a_{3^\prime(0)}T(a_{4^\prime})_{(0)}a_{6^\prime}\\
&=&T(a_{1^\prime})_{(-1)1^\prime}a_{2^\prime}\otimes T(a_{1^\prime})_{(-1)2^\prime}a_{3^\prime(-1)}T(a_{4^\prime})_{(-1)}a_{5^\prime}\otimes T(a_{1^\prime})_{(0)}a_{3^\prime(0)}T(a_{4^\prime})_{(0)}a_{6^\prime}\\
&=&T(a_{1^\prime})_{(-1)1^\prime}a_{2^\prime}\otimes T(a_{1^\prime})_{(-1)2^\prime}a_{3^\prime}\otimes T(a_{1^\prime})_{(0)}a_{4^\prime}\\
&=&a_{[0]1^\prime}\otimes a_{[0]2^\prime}\otimes a_{[1]},
\end{eqnarray*}
so, Eq.$(HM3)$ holds. In addition, $a_{(-1)}\v(a_{(0)})=\v(a)1_A$ by Lemma 2.8, and $$\v(a_{[0]})a_{[1]}=\v(T(a_{1^\prime})_{(-1)}a_{2^\prime})T(a_{1^\prime})_{(0)}a_{3^\prime}=T(a_{1^\prime})a_{2^\prime}=\v(a)1_A,
$$
that is, Eq.$(HM1)$ holds, which completes the proof.\hfill $\square$

\vspace{2mm}

\noindent{\bf Proposition 3.4}\ \
Let $(A,\Delta^\prime)$ be a commutative Hopf algebra with antipode $T$. Assume that $(A,A)$ is a Hopf matched pair with coactions $\rho$ and $\varphi$, such that  for all $a \in A$,
$$a_{1^\prime}\otimes a_{2^\prime}=a_{1^\prime(-1)}a_{2^\prime[0]}\otimes a_{1^\prime(0)}a_{2^\prime[1]}.\eqno(3.1)
$$

Then, $(A,\Delta,\Delta^\prime)$ is a commutative Hopf brace with
$$\Delta(a)\equiv a_1\otimes a_2=a_{1^\prime}T(a_{2^\prime(-1)})\otimes a_{2^\prime(0)},\eqno{(3.2)}$$
$$S(a)=a_{(-1)}T(a_{(0)}),\eqno{(3.3)}
$$
for all $a \in A$.

\noindent{\bf Proof.}
Since $\rho$ and $\Delta^\prime$ are algebra maps and $A$ is commutative, it is easy to see that $\D$ is an algebra map.

Further, for any $a \in A$, we have
\begin{eqnarray*}
a_1a_{2(-1)}\otimes a_{2(0)}&\stackrel{(3.2)}{=}&a_{1^\prime}T(a_{2^\prime(-1)})a_{2^\prime(0)(-1)}\otimes a_{2^\prime(0)(0)}\\
&=&a_{1^\prime}T(a_{2^\prime(-1)1^\prime})a_{2^\prime(-1)2^\prime}\otimes a_{2^\prime(0)}\\
&=&a_{1^\prime}\varepsilon(a_{2^\prime(-1)})\otimes a_{2^\prime(0)}\\
&=&a_{1^\prime}\otimes a_{2^\prime}.
\end{eqnarray*}

Thus, we get
$$a_{1^\prime}\otimes a_{2^\prime}=a_1a_{2(-1)}\otimes a_{2(0)}.\eqno(3.4)$$

According to Eq.(3.4), we have
\begin{eqnarray*}
(id\otimes \Delta)\rho(a)&=&a_{(-1)}\otimes a_{(0)1}\otimes a_{(0)2}\\
&\stackrel{(3.2)}{=}&a_{(-1)}\otimes a_{(0)1^\prime}T(a_{(0)2^\prime(-1)})\otimes a_{(0)2^\prime(0)}\\
&\stackrel{(HMC2)}{=}&a_{1^\prime(-1)}a_{2^\prime(-1)[0]}\otimes a_{1^\prime(0)}a_{2^\prime(-1)[1]}T(a_{2^\prime(0)(-1)})\otimes a_{2^\prime(0)(0)}\\
&\stackrel{(3.4)}{=}&a_{1(-1)}a_{2(-1)(-1)}a_{2(0)(-1)[0]}\otimes a_{1(0)}a_{2(-1)(0)}a_{2(0)(-1)[1]}T(a_{2(0)(0)(-1)})\otimes a_{2(0)(0)(0)}\\
&=&a_{1(-1)}a_{2(-1)1^\prime(-1)}a_{2(-1)2^\prime[0]}\otimes a_{1(0)}a_{2(-1)1^\prime(0)}a_{2(-1)2^\prime[1]}T(a_{2(0)(-1)})\otimes a_{2(0)(0)}\\
&\stackrel{(HMC1)}{=}&a_{1(-1)}a_{2(-1)1^\prime}\otimes a_{1(0)}a_{2(-1)2^\prime}T(a_{2(0)(-1)})\otimes a_{2(0)(0)}\\
&=&a_{1(-1)}a_{2(-1)}\otimes a_{1(0)}a_{2(0)(-1)}T(a_{2(0)(0)(-1)})\otimes a_{2(0)(0)(0)}\\
&=&a_{1(-1)}a_{2(-1)}\otimes a_{1(0)}a_{2(0)(-1)1^\prime}T(a_{2(0)(-1)2^\prime})\otimes a_{2(0)(0)}\\
&=&a_{1(-1)}a_{2(-1)}\otimes a_{1(0)}\varepsilon(a_{2(0)(-1)})\otimes a_{2(0)(0)}\\
&=&a_{1(-1)}a_{2(-1)}\otimes a_{1(0)}\otimes a_{2(0)}.
\end{eqnarray*}

Thus, we get
$$a_{(-1)}\o a_{(0)1}\o a_{(0)2}=a_{1(-1)}a_{2(-1)}\otimes a_{1(0)}\otimes a_{2(0)}.\eqno(3.5)
$$

In what follows, according to Eq.(3.4) and Eq.(3.5), we prove that $(A,\Delta,\varepsilon,S)$ is a Hopf algebra.

As a matter of fact, we have
\begin{eqnarray*}
(id\o \D)\D(a)&\stackrel{(3.2)}{=}&a_{1^\prime}T(a_{2^\prime(-1)})\otimes a_{2^\prime(0)1}\otimes a_{2^\prime(0)2}\\
&\stackrel{(3.5)}{=}&a_{1^\prime}T(a_{2^\prime1(-1)}a_{2^\prime2(-1)})\otimes a_{2^\prime1(0)}\otimes a_{2^\prime2(0)}\\
&=&a_{1^\prime}T(a_{2^\prime2(-1)})T(a_{2^\prime1(-1)})\otimes a_{2^\prime1(0)}\otimes a_{2^\prime2(0)}\\
&\stackrel{(3.2)}{=}&a_{1^\prime}T(a_{3^\prime(0)(-1)})T((a_{2^\prime}T(a_{3^\prime(-1)}))_{(-1)})\otimes (a_{2^\prime}T(a_{3^\prime(-1)}))_{(0)}\otimes a_{3^\prime(0)(0)}\\
&\stackrel{(3.4)}{=}&a_{1^\prime1}a_{1^\prime2(-1)}T(a_{2^\prime(0)(-1)})T((a_{1^\prime2(0)}T(a_{2^\prime(-1)}))_{(-1)})\\
&&\otimes (a_{1^\prime2(0)}T(a_{2^\prime(-1)}))_{(0)}\otimes a_{2^\prime(0)(0)}\\
&=&a_{1^\prime1}a_{1^\prime2(-1)}T(a_{2^\prime(-1)2^\prime})T(T(a_{2^\prime(-1)1^\prime})_{(-1)})T(a_{1^\prime2(0)(-1)})\\
&&\otimes a_{1^\prime2(0)(0)}T(a_{2^\prime(-1)1^\prime})_{(0)}\otimes a_{2^\prime(0)}\\
&=&a_{1^\prime1}a_{1^\prime2(-1)1^\prime}T(a_{2^\prime(-1)2^\prime})T(T(a_{2^\prime(-1)1^\prime})_{(-1)})T(a_{1^\prime2(-1)2^\prime})\\
&&\otimes a_{1^\prime2(0)}T(a_{2^\prime(-1)1^\prime})_{(0)}\otimes a_{2^\prime(0)}\\
\end{eqnarray*}
\begin{eqnarray*}
&=&a_{1^\prime1}a_{1^\prime2(-1)1^\prime}T(a_{1^\prime2(-1)2^\prime})T(a_{2^\prime(-1)2^\prime})T(T(a_{2^\prime(-1)1^\prime})_{(-1)})\\
&&\otimes a_{1^\prime2(0)}T(a_{2^\prime(-1)1^\prime})_{(0)}\otimes a_{2^\prime(0)}\\
&=&a_{1^\prime1}T(a_{2^\prime(-1)2^\prime})T(T(a_{2^\prime(-1)1^\prime})_{(-1)})\otimes a_{1^\prime2}T(a_{2^\prime(-1)1^\prime})_{(0)}\otimes a_{2^\prime(0)}\\
&=&a_{1^\prime1}T(a_{2^\prime(-1)})_{1^\prime}T(T(a_{2^\prime(-1)})_{2^\prime(-1)})\otimes a_{1^\prime2}T(a_{2^\prime(-1)})_{2^\prime(0)}\otimes a_{2^\prime(0)}\\
&\stackrel{(3.2)}{=}&a_{1^\prime1}T(a_{2^\prime(-1)})_{1}\otimes a_{1^\prime2}T(a_{2^\prime(-1)})_2\otimes a_{2^\prime(0)}\\
&=&(a_{1^\prime}T(a_{2^\prime(-1)}))_{1}\otimes (a_{1^\prime}T(a_{2^\prime(-1)}))_2\otimes a_{2^\prime(0)}\\
&=&a_{11}\otimes a_{12}\otimes a_2\\&=&(\D\o id)\D(a),\\
(id\otimes \varepsilon)\Delta(a)&=&a_{1^\prime}T(a_{2^\prime(-1)})\varepsilon(a_{2^\prime(0)})\\
&\stackrel{(HMC2)}=&a_{1^\prime}\varepsilon(a_{2^\prime})=a=(\varepsilon\otimes id)\Delta(a),
\end{eqnarray*}
so, $(A,\Delta,\varepsilon)$ is a coalgebra. Moreover, we have
\begin{eqnarray*}
a_1S(a_2)&\stackrel{(3.3)}{=}&a_1a_{2(-1)}T(a_{2(0)})\\
&\stackrel{(3.2)}{=}&a_{1^\prime}T(a_{2^\prime(-1)})a_{2^\prime(0)(-1)}T(a_{2^\prime(0)(0)})\\
&=&a_{1^\prime}T(a_{2^\prime(-1)1^\prime})a_{2^\prime(-1)2^\prime}T(a_{2^\prime(0)})\\
&=&a_{1^\prime}\varepsilon(a_{2^\prime(-1)})T(a_{2^\prime(0)})\\
&=&a_{1^\prime}T(a_{2^\prime})=\varepsilon(a)1
\end{eqnarray*}
for all $a \in A,$ so, by the commutativity of $A$, we easily know that $(A,m,\varepsilon,\Delta,1,S)$ is a Hopf algebra.

Again by Eq.(3.4), we get
\begin{eqnarray*}&&~~~~S(a_1)a_{21^\prime}\otimes a_{22^\prime}=S(a_1)a_2a_{3(-1)}\otimes a_{3(0)}=a_{(-1)}\otimes a_{(0)}=\rho(a)~~~~~~~~(3.6)
\end{eqnarray*} for each $a \in A.$

Finally, we have only to prove the compatible condition of Hopf braces is satisfied:
\begin{eqnarray*}
a_{11^\prime}S(a_2)a_{31^\prime}\otimes a_{12^\prime}\otimes a_{32^\prime}&\stackrel{(3.6)}{=}&a_{11^\prime}a_{2(-1)}\otimes a_{12^\prime}\otimes a_{2(0)}\\
&\stackrel{(3.2)}{=}&a_{1^\prime}T(a_{3^\prime(-1)})_{1^\prime}a_{3^\prime(0)(-1)}\otimes a_{2^\prime}T(a_{3^\prime(-1)})_{2^\prime}\otimes a_{3^\prime(0)(0)}\\
&=&a_{1^\prime}T(a_{3^\prime(-1)1^\prime})_{1^\prime}a_{3^\prime(-1)2^\prime}\otimes a_{2^\prime}T(a_{3^\prime(-1)1^\prime})_{2^\prime}\otimes a_{3^\prime(0)}\\
&=&a_{1^\prime}T(a_{3^\prime(-1)2^\prime})a_{3^\prime(-1)3^\prime}\otimes a_{2^\prime}T(a_{3^\prime(-1)1^\prime})\otimes a_{3^\prime(0)}\\
&=&a_{1^\prime}\varepsilon(a_{3^\prime(-1)2^\prime})\otimes a_{2^\prime}T(a_{3^\prime(-1)1^\prime})\otimes a_{3^\prime(0)}\\
&=&a_{1^\prime}\otimes a_{2^\prime}T(a_{3^\prime(-1)})\otimes a_{3^\prime(0)}\\
&=&a_{1^\prime}\otimes a_{2^\prime1^\prime}T(a_{2^\prime2^\prime(-1)})\otimes a_{2^\prime2^\prime(0)}\\
&=&a_{1^\prime}\otimes a_{2^\prime1}\otimes a_{2^\prime2}
\end{eqnarray*} for each $a \in A.$
 \hfill $\square$

\vspace{2mm}

 Let $(A,\Delta)$ be a commutative Hopf algebra with antipode $S$. Let $\mathcal{M}(A)$ be the category with objects Hopf matched pairs $(A,A)$ such that the condition (3.1) is satisfied, and all morphisms Hopf algebra homomorphisms $f:A\rightarrow A$ such that $\rho f(a)=(f\otimes f)\rho(a),~\varphi f(a)=(f\otimes f)\varphi(a),$ for all $a \in A$.

\vspace{2mm}

\noindent{\bf Theorem 3.5}\ \
Let $(A,\Delta)$ be a commutative Hopf algebra with antipode $S$. Then, the category $\mathcal{CB}(A)$ of Hopf braces is equivalent to the category $\mathcal{M}(A)$ of Hopf matched pairs.

\noindent{\bf Proof.}\ \ We have two functors as follows:
\begin{eqnarray*}
F:\mathcal{CB}(A)&\rightarrow&\mathcal{M}(A), \ \ F((A,\Delta,\Delta^\prime))=(A,A),\\
F(f)&=&f,
\end{eqnarray*}
  where $(A,A)$ is the Hopf matched pair as in  Proposition 3.3.
  \begin{eqnarray*}
F:\mathcal{M}(A)&\rightarrow&\mathcal{CB}(A), \ \ G((A,A))=(A,\Delta,\Delta^\prime),\\
G(f)&=&f,
\end{eqnarray*}
  where $(A,\Delta,\Delta^\prime)$ is a Hopf brace as in Proposition 3.4.

By a direct calculation, we can show that $\mathcal{CB}(A)$ is equivalent to $\mathcal{M}(A)$. \hfill $\square$

\vspace{3mm}

  \begin{center}
{\bf \S4\quad Hopf brace on bicrossed coproduct}
\end{center}
\vspace{3mm}

In this section, we mainly construct Hopf braces on bicrossed coproducts. More precisely, we give a sufficient and necessary condition for a given bicrossed coproduct $A\bowtie H$ to be a Hopf brace if $A$ or $H$ is a Hopf brace, and show that the dual of Drinfel'd double $D(H)$ of a finite dimensional cocommutative Hopf algebra $H$ is a Hopf brace.

 Assume that $(A,H,\rho,\varphi)$ is a Hopf matched pair in Definition 3.1, and $A\o H$ a $k$-module with the tensor algebra structure. Define a comultiplication on  $A\o H$  as follows: for all $a\in A, h\in H$,
$$
\widetilde{\Delta}_{A\otimes H}(a\otimes h)=a_1\otimes a_{2(-1)}h_{1[0]}\otimes a_{2(0)}h_{1[1]}\otimes h_2.
$$

Then, by \cite{H. X. Chen}, $A\o H$ is a Hopf algebra, whose antipode is given by
$$
\widetilde{S}(a\o h)=S_A(h_{[1]})S_A(a_{(0)})\o S_H(h_{[0]})S_H(a_{(-1)}).
$$

In what follows, we call the Hopf algebra a bicrossed coproduct of $A$ and $H$, and denote it by $A\bowtie H$, whose  comultiplication is denoted by $\widetilde{\Delta}.$

\vspace{2mm}

\noindent {\bf Proposition 4.1}\ \ Let $(A,\Delta,\Delta^\prime)$ be a Hopf brace, and $H$ a commutative cocommutative Hopf algebra. If $(A,\rho)$ is a left $H$-comodule coalgebra, and $(A_{\Delta^\prime},H,\rho,\varphi)$ a Hopf matched pair,
where $A_{\Delta^\prime}$ denotes the Hopf algebra with comultiplication $\Delta^\prime$ of $A$. Then, $(A_{\Delta^\prime}\bowtie H, \widehat{\Delta}, \widetilde{\Delta})$ is a Hopf brace, if and only if
$$
h_{[0]}\otimes h_{[1]1}\otimes h_{[1]2}=h_{1[0]}S(h_2)h_{3[0]}\otimes h_{1[1]}\otimes h_{3[1]},\eqno(4.1)
$$
for all $h\in H$, where $\widehat{\Delta}$ denotes the comultiplication of the usual tensor coalgebra of $A\otimes H$.

\noindent{\bf Proof.}\ \ $``\Longrightarrow"$  Let the condition (4.1) holds. In what follows, we have only to prove Eq.(2.1) is satisfied in order to show $(A_{\Delta^\prime}\bowtie H, \widehat{\Delta}, \widetilde{\Delta})$ to be a Hopf brace.

For the sake of convenience, we denote $\widehat{\Delta}(a\otimes h)$ by $(a\otimes h)_{\widehat{1}}\otimes (a\otimes h)_{\widehat{2}}$ and $\widetilde{\Delta}(a\otimes h)$ by $(a\otimes h)_{\widetilde{1}}\otimes (a\otimes h)_{\widetilde{2}}$.

Thus, for all $a\in A $ and $h\in H$, we have
\begin{eqnarray*}
&&(a\otimes h)_{\widetilde 1}\otimes (a\otimes h)_{{\widetilde 2}{\widehat 1}}\otimes (a\otimes h)_{{\widetilde 2}{\widehat 2}}\\
&=&a_{1^\prime}\otimes a_{2^\prime(-1)}h_{1[0]}\otimes a_{2^\prime(0)1}h_{1[1]1}\otimes h_2 \otimes a_{2^\prime(0)2}h_{1[1]2}\otimes h_3\\
&=&a_{1^\prime}\otimes a_{2^\prime1(-1)}a_{2^\prime2(-1)}h_{1[0]}\otimes a_{2^\prime1(0)}h_{1[1]1}\otimes h_2\otimes a_{2^\prime2(0)}h_{1[1]2}\otimes h_3\\
&\stackrel{(2.1)}{=}&a_{11^\prime}S(a_2)a_{31^\prime}\otimes a_{12^\prime(-1)}a_{32^\prime(-1)}h_{1[0]}\otimes a_{12^\prime(0)}h_{1[1]1}\otimes h_2\otimes a_{32^\prime(0)}h_{1[1]2}\otimes h_3\\
&\stackrel{(4.1)}{=}&a_{11^\prime}S(a_2)a_{31^\prime}\otimes a_{12^\prime(-1)}a_{32^\prime(-1)}h_{1[0]}S(h_2)h_{3[0]}\otimes a_{12^\prime(0)}h_{1[1]}\otimes h_4\otimes a_{32^\prime(0)}h_{3[1]}\otimes h_5\\
&=&a_{11^\prime}S(a_2)a_{31^\prime}\otimes a_{12^\prime(-1)}h_{1[0]}S(h_2)a_{32^\prime(-1)}h_{3[0]}\otimes a_{12^\prime(0)}h_{1[1]}\otimes h_4\otimes a_{32^\prime(0)}h_{3[1]}\otimes h_5\\
&=&a_{11^\prime}S(a_2)a_{31^\prime}\otimes a_{12^\prime(-1)}h_{1[0]}S(h_3)a_{32^\prime(-1)}h_{4[0]}\otimes a_{12^\prime(0)}h_{1[1]}\otimes h_2\otimes a_{32^\prime(0)}h_{4[1]}\otimes h_5\\
&=&(a\otimes h)_{\widehat 1 \widetilde 1}S((a\otimes h)_{\widehat 2})(a\otimes h)_{\widehat 3 \widetilde 1}\otimes (a\otimes h)_{\widehat 1 \widetilde 2}\otimes (a\otimes h)_{\widehat 3 \widetilde2}.
\end{eqnarray*}

$``\Longleftarrow"$ If $(A_{\Delta^\prime}\bowtie H, \widehat{\Delta}, \widetilde{\Delta})$ is a Hopf brace, then, by the above proof, we have
\begin{eqnarray*}
&&a_{11^\prime}S(a_2)a_{31^\prime}\otimes a_{12^\prime(-1)}a_{32^\prime(-1)}h_{1[0]}\otimes a_{12^\prime(0)}h_{1[1]1}\otimes h_2\otimes a_{32^\prime(0)}h_{1[1]2}\otimes h_3\\
&=&a_{11^\prime}S(a_2)a_{31^\prime}\otimes a_{12^\prime(-1)}h_{1[0]}S(h_2)a_{32^\prime(-1)}h_{3[0]}\otimes a_{12^\prime(0)}h_{1[1]}\otimes h_4\otimes a_{32^\prime(0)}h_{3[1]}\otimes h_5
\end{eqnarray*}
for any $a\in A, h\in H$.

If taking $a=1$, then, the above equation translates into the following equation
\begin{eqnarray*}
&&1\otimes h_{1[0]}\otimes h_{1[1]1}\otimes h_2\otimes h_{1[1]2}\otimes h_3\\
&=&1 \otimes h_{1[0]}S(h_2)h_{3[0]}\otimes h_{1[1]}\otimes h_4\otimes h_{3[1]}\otimes h_5.
\end{eqnarray*}

By applying $\v\o id\o id \o \v\o id\o \v$ to both sides of the above equation, we get Eq(4.1).
\hfill $\square$

\vspace{2mm}

\noindent{\bf Remark 4.2}\ \ (1) Let $(A,\Delta,\Delta^\prime)$ be a Hopf brace, and $H$ a commutative cocommutative Hopf algebra. Suppose that the right $A_{\Delta^\prime}$-comodule action of $H$ is trivial. Then, by Definition 3.1, $(A_{\Delta^\prime},H,\rho,\varphi)$ is a Hopf matched pair if $(A_{\Delta^\prime},\rho)$ is a left $H$-comodule bialgebra.

 It is obvious that the condition (4.1) holds. So, according to Proposition 4.1, $(A_{\Delta^\prime}\bowtie H, \widehat{\Delta}, \widetilde{\Delta})$ is a Hopf brace if $(A,\rho)$ is a left $H$-comudule coalgebra, where the comultiplication $\widetilde{\Delta}$ is given by
$$
\widetilde{\Delta}(a\otimes h)=a_{1^\prime}\otimes a_{2^\prime(-1)}h_1\otimes a_{2^\prime(0)}\otimes h_2.
$$

In this case, the comultiplication $\widetilde{\Delta}$ of the bicrossed coproduct $A_{\Delta^\prime}\bowtie H$ is actually the comultiplication of the usual smash coproduct on $A_{\Delta^\prime}\otimes H$.

(2) Suppose that $A$ is a Hopf algebra with comultiplication $\Delta$. Then, $(A,\Delta,\Delta)$ is a Hopf brace. If $H$ is a commutative cocommutative Hopf algebra, and $(A,\rho)$ is a left $H$-comudule bialgebra. Then, by the above remark, we know that the smash coproduct $(A\times H, \widehat{\Delta}, \widetilde{\Delta})$ is a Hopf brace.

\vspace{2mm}

\noindent{\bf Proposition 4.3}\ \ Let $A$ be a Hopf algebra, and $(H,\Delta,\Delta^\prime)$ a commutative Hopf brace. If $(A,H,\rho,\varphi)$ is a Hopf matched pair, and $(A,\rho^\prime)$ a left $H_{\Delta^\prime}$-comodule bialgebra (whose comodule structure is given by $\rho^\prime(a)=a_{(-1)^\prime}\o a_{(0)^\prime}$ for $a\in A$). Then, $(A\bowtie H, \widetilde{\Delta},\bar{\D})$ is a Hopf brace, if and only if for all $a\in A,h\in H$,
\begin{eqnarray*}
&&a_{(-1)^\prime}\otimes a_{(0)^\prime(-1)}\otimes a_{(0)^\prime(0)}=a_{(-1)11^\prime}S(a_{(-1)2})a_{(0)(-1)^\prime}\otimes a_{(-1)12^\prime}\otimes a_{(0)(0)^\prime},~~~~~~~~~~~~~~~~~(4.2)\\
&&h_{1^\prime}\otimes h_{2^\prime[0]}\otimes h_{2^\prime[1]}=h_{1[0]11^\prime}S(h_{1[0]2})h_{1[1](-1)^\prime}h_2\otimes h_{1[0]12^\prime}\otimes h_{1[1](0)^\prime},~~~~~~~~~~~~~~~~~~~~~~~~~~(4.3)
\end{eqnarray*}
where $\bar{\Delta}$ denotes the comultiplication of smash coproduct on $A\otimes H_{\D'}$, that is,
$$\begin{array}{rllr}
\bar{\Delta}(a\o h)&\equiv (a\otimes h)_{\bar{1}}\otimes (a\otimes h)_{\bar{2}}\\
&=a_1\o a_{2(-1)'}h_{1'}\o a_{2(0)'}\o h_{2'}.
\end{array}$$

\noindent{\bf Proof.}\ \ $``\Longrightarrow"$ Since $(A,\rho^\prime)$ is a left $H_{\Delta^\prime}$-comodule bialgebra, the smash coproduct with the comultiplication $\bar{\Delta}$ on the tensor product algebra $A\otimes H_{\D'}$  is also a Hopf algebra. Therefore, we have only  to prove that Eq.(2.1) holds.

As a matter of fact, for all $a \in A,h \in H$, we have
\begin{eqnarray*}
&&(a\otimes h)_{\bar1}\otimes (a\otimes h)_{\bar2 \widetilde1}\otimes (a\otimes h)_{\bar2\widetilde2}\\
&=&a_1\otimes \underline{a_{2(-1)^\prime}}h_{1^\prime}\otimes \underline{a_{2(0)^\prime1}}\otimes \underline{a_{2(0)^\prime2(-1)}}h_{2^\prime1[0]}\otimes a_{2(0)^\prime2(0)}h_{2^\prime1[1]}\otimes h_{2^\prime2}\\
&=&a_1\otimes \underline{a_{2(-1)^\prime}a_{3(-1)^\prime}}h_{1^\prime}\otimes \underline{a_{2(0)^\prime}}\otimes \underline{a_{3(0)^\prime(-1)}}h_{2^\prime1[0]}\otimes a_{3(0)^\prime(0)}h_{2^\prime1[1]}\otimes h_{2^\prime2}\\
&\stackrel{(4.2)}{=}&a_1\otimes a_{2(-1)^\prime}a_{3(-1)11^\prime}S(a_{3(-1)2})a_{3(0)(-1)^\prime}h_{1^\prime}\otimes a_{2(0)^\prime}\otimes a_{3(-1)12^\prime}h_{2^\prime1[0]}\\
&&\otimes a_{3(0)(0)^\prime}h_{2^\prime1[1]}\otimes h_{2^\prime2}\\
&\stackrel{(2.1)}{=}&a_1\otimes a_{2(-1)^\prime}a_{3(-1)11^\prime}S(a_{3(-1)2})a_{3(0)(-1)^\prime}h_{11^\prime}S(h_2)h_{31^\prime}\otimes a_{2(0)^\prime}\otimes a_{3(-1)12^\prime}h_{12^\prime[0]}\\
&&\otimes a_{3(0)(0)^\prime}h_{12^\prime[1]}\otimes h_{32^\prime}\\
&\stackrel{(4.3)}{=}&a_1\otimes a_{2(-1)^\prime}a_{3(-1)11^\prime}S(a_{3(-1)2})a_{3(0)(-1)^\prime}h_{1[0]11^\prime}S(h_{1[0]2})h_{1[1](-1)^\prime}h_2S(h_3)h_{41^\prime}\\
&&\otimes a_{2(0)^\prime}\otimes a_{3(-1)12^\prime}h_{1[0]12^\prime}\otimes a_{3(0)(0)^\prime}h_{1[1](0)^\prime}\otimes h_{42^\prime}\\
&=&a_1\otimes a_{2(-1)^\prime}a_{3(-1)11^\prime}S(a_{3(-1)2})a_{3(0)(-1)^\prime}h_{1[0]11^\prime}S(h_{1[0]2})h_{1[1](-1)^\prime}h_{21^\prime}\\
&&\otimes a_{2(0)^\prime}\otimes a_{3(-1)12^\prime}h_{1[0]12^\prime}\otimes a_{3(0)(0)^\prime}h_{1[1](0)^\prime}\otimes h_{22^\prime}\\
&=&a_1\otimes a_{2(-1)^\prime}a_{3(-1)11^\prime}h_{1[0]11^\prime}S(h_{1[0]2})S(a_{3(-1)2})a_{3(0)(-1)^\prime}h_{1[1](-1)^\prime}h_{21^\prime}\\&&
\otimes a_{2(0)^\prime}\otimes a_{3(-1)12^\prime}h_{1[0]12^\prime}\otimes a_{3(0)(0)^\prime}h_{1[1](0)^\prime}\otimes h_{22^\prime}\\
&=&a_1S(h_{1[1]1})h_{1[1]2}\otimes a_{2(-1)^\prime}a_{3(-1)11^\prime}h_{1[0]11^\prime}S(h_{1[0]2})S(a_{3(-1)2})a_{3(0)(-1)^\prime}h_{1[1]3(-1)^\prime}
h_{21^\prime}\\&&
\otimes a_{2(0)^\prime}\otimes a_{3(-1)12^\prime}h_{1[0]12^\prime}\otimes a_{3(0)(0)^\prime}h_{1[1]3(0)^\prime}\otimes h_{22^\prime}\\&=&a_1S(h_{1[0][1]})h_{1[1]1}\otimes a_{2(-1)^\prime}a_{3(-1)11^\prime}h_{1[0][0]11^\prime}S(h_{1[0][0]2})S(a_{3(-1)2})a_{3(0)(-1)^\prime}h_{1[1]2(-1)^\prime}
\\&& h_{21^\prime}\otimes a_{2(0)^\prime}\otimes a_{3(-1)12^\prime}h_{1[0][0]12^\prime}\otimes a_{3(0)(0)^\prime}h_{1[1]2(0)^\prime}\otimes h_{22^\prime}\\
&=&a_1S(h_{1[0][1]})S(a_{3(0)1})a_{3(0)2}h_{1[1]1}\otimes a_{2(-1)^\prime}a_{3(-1)11^\prime}h_{1[0][0]11^\prime}S(h_{1[0][0]2})S(a_{3(-1)2})
\\&&a_{3(0)3(-1)^\prime}h_{1[1]2(-1)^\prime}
h_{21^\prime}\otimes a_{2(0)^\prime}\otimes a_{3(-1)12^\prime}h_{1[0][0]12^\prime}\otimes a_{3(0)3(0)^\prime}h_{1[1]2(0)^\prime}\otimes h_{22^\prime}\\
&\stackrel{(3.2)}{=}&a_1S(h_{1[0][1]})S(a_{3(0)}a_{4(-1)[1]})a_{4(0)1}h_{1[1]1}\otimes a_{2(-1)^\prime}a_{3(-1)11^\prime}a_{4(-1)[0]11^\prime}h_{1[0][0]11^\prime}\\
&&S(h_{1[0][0]2}) S(a_{3(-1)2}a_{4(-1)[0]2})a_{4(0)2(-1)^\prime}h_{1[1]2(-1)^\prime}
h_{21^\prime}\otimes a_{2(0)^\prime}\otimes a_{3(-1)12^\prime}a_{4(-1)[0]12^\prime}\\
&&h_{1[0][0]12^\prime}\otimes a_{4(0)2(0)^\prime}h_{1[1]2(0)^\prime}\otimes h_{22^\prime}\\
&=&a_1S(h_{1[0][1]})S(a_{3(0)(0)}a_{4(-1)[1]})a_{4(0)1}h_{1[1]1}\otimes a_{2(-1)^\prime}a_{3(-1)1^\prime}a_{4(-1)[0]11^\prime}h_{1[0][0]11^\prime}\\
&&S(h_{1[0][0]2})S(a_{3(0)(-1)}a_{4(-1)[0]2})a_{4(0)2(-1)^\prime}h_{1[1]2(-1)^\prime}
h_{21^\prime}\otimes a_{2(0)^\prime}\otimes a_{3(-1)2^\prime}a_{4(-1)[0]12^\prime}\\
&&h_{1[0][0]12^\prime}\otimes a_{4(0)2(0)^\prime}h_{1[1]2(0)^\prime}\otimes h_{22^\prime}\\
\end{eqnarray*}
\begin{eqnarray*}
&\stackrel{(3.3)}{=}&a_1S(h_{1[0][1]})S(a_{3(0)(0)}a_{4(-1)1[1](0)}a_{4(-1)2[1]})a_{4(0)1}h_{1[1]1}\otimes
a_{2(-1)^\prime}a_{3(-1)1^\prime}a_{4(-1)1[0]1^\prime}\\
&&h_{1[0][0]11^\prime}S(h_{1[0][0]2})S(a_{3(0)(-1)}a_{4(-1)[1](-1)}a_{4(-1)2[0]})a_{4(0)2(-1)^\prime}h_{1[1]2(-1)^\prime}
h_{21^\prime}\\
&&\otimes a_{2(0)^\prime}\otimes a_{3(-1)2^\prime}a_{4(-1)1[0]2^\prime}h_{1[0][0]12^\prime}\otimes a_{4(0)2(0)^\prime}h_{1[1]2(0)^\prime}\otimes h_{22^\prime}\\
&\stackrel{(3.3)}{=}&a_1S(h_{1[0]1[1](0)}h_{1[0]2[1]})S(a_{3(0)(0)}a_{4(-1)1[1](0)}a_{4(-1)2[1]})a_{4(0)1}h_{1[1]1}\otimes a_{2(-1)^\prime}a_{3(-1)1^\prime}\\&&a_{4(-1)1[0]1^\prime}h_{1[0]1[0]1^\prime}S(h_{1[0]1[1](-1)}h_{1[0]2[0]})S(a_{3(0)(-1)}a_{4(-1)[1](-1)}a_{4(-1)2[0]})a_{4(0)2(-1)^\prime}
\\
&&h_{1[1]2(-1)^\prime}
h_{21^\prime}\otimes a_{2(0)^\prime}\otimes a_{3(-1)2^\prime}a_{4(-1)1[0]2^\prime}h_{1[0]1[0]2^\prime}\otimes a_{4(0)2(0)^\prime}h_{1[1]2(0)^\prime}\otimes h_{22^\prime}\\
&=&(a\otimes h)_{\widetilde1\bar1}S((a\otimes h)_{\widetilde2})(a\otimes h)_{\widetilde3\bar1}\otimes (a\otimes h)_{\widetilde1\bar2}\otimes (a\otimes h)_{\widetilde3\bar2},
\end{eqnarray*}
which is complete.

$``\Longleftarrow"$ If $(A\bowtie H, \widetilde{\Delta},\bar{\D})$ is a Hopf brace, then, by the above proof, we get
\begin{eqnarray*}
&&a_1\otimes a_{2(-1)^\prime}a_{3(-1)^\prime}h_{1^\prime}\otimes a_{2(0)^\prime}\otimes a_{3(0)^\prime(-1)}h_{2^\prime1[0]}\otimes a_{3(0)^\prime(0)}h_{2^\prime1[1]}\otimes h_{2^\prime2}\\
&=&a_1\otimes a_{2(-1)^\prime}a_{3(-1)11^\prime}S(a_{3(-1)2})a_{3(0)(-1)^\prime}h_{1[0]11^\prime}S(h_{1[0]2})h_{1[1](-1)^\prime}h_{21^\prime}\\
&&\otimes a_{2(0)^\prime}\otimes a_{3(-1)12^\prime}h_{1[0]12^\prime}\otimes a_{3(0)(0)^\prime}h_{1[1](0)^\prime}\otimes h_{22^\prime}.
\end{eqnarray*}

By setting $a=1$ and then applying $\v\o id\o \v\o id\o id\o \v$ to both sides of the equation, we get Eq.(4.3). So, again by the above proof, we have
\begin{eqnarray*}
&&a_1\otimes a_{2(-1)^\prime}a_{3(-1)^\prime}h_{1^\prime}\otimes a_{2(0)^\prime}\otimes a_{3(0)^\prime(-1)}h_{2^\prime1[0]}\otimes a_{3(0)^\prime(0)}h_{2^\prime1[1]}\otimes h_{2^\prime2}\\
&=&a_1\otimes a_{2(-1)^\prime}a_{3(-1)11^\prime}S(a_{3(-1)2})a_{3(0)(-1)^\prime}h_{11^\prime}S(h_2)h_{31^\prime}\otimes a_{2(0)^\prime}\otimes a_{3(-1)12^\prime}h_{12^\prime[0]}\\
&&\otimes a_{3(0)(0)^\prime}h_{12^\prime[1]}\otimes h_{32^\prime}.
\end{eqnarray*}

If Eq.(4.3) holds, we get
\begin{eqnarray*}
&&a_1\otimes a_{2(-1)^\prime}a_{3(-1)11^\prime}S(a_{3(-1)2})a_{3(0)(-1)^\prime}h_{1[0]11^\prime}S(h_{1[0]2})h_{1[1](-1)^\prime}h_{21^\prime}\\
&&\otimes a_{2(0)^\prime}\otimes a_{3(-1)12^\prime}h_{1[0]12^\prime}\otimes a_{3(0)(0)^\prime}h_{1[1](0)^\prime}\otimes h_{22^\prime}\\
&=&a_1\otimes a_{2(-1)^\prime}a_{3(-1)11^\prime}S(a_{3(-1)2})a_{3(0)(-1)^\prime}h_{1[0]11^\prime}S(h_{1[0]2})h_{1[1](-1)^\prime}h_2S(h_3)h_{41^\prime}\\
&&\otimes a_{2(0)^\prime}\otimes a_{3(-1)12^\prime}h_{1[0]12^\prime}\otimes a_{3(0)(0)^\prime}h_{1[1](0)^\prime}\otimes h_{42^\prime}\\
&\stackrel{(4.3)}{=}&a_1\otimes a_{2(-1)^\prime}a_{3(-1)11^\prime}S(a_{3(-1)2})a_{3(0)(-1)^\prime}h_{11^\prime}S(h_2)h_{31^\prime}\otimes a_{2(0)^\prime}\otimes a_{3(-1)12^\prime}h_{12^\prime[0]}\\
&&\otimes a_{3(0)(0)^\prime}h_{12^\prime[1]}\otimes h_{32^\prime}\\
&\stackrel{(2.1)}{=}&a_1\otimes a_{2(-1)^\prime}a_{3(-1)11^\prime}S(a_{3(-1)2})a_{3(0)(-1)^\prime}h_{1^\prime}\otimes a_{2(0)^\prime}\otimes a_{3(-1)12^\prime}h_{2^\prime1[0]}\\
&&\otimes a_{3(0)(0)^\prime}h_{12^\prime[1]}\otimes h_{2^\prime2},
\end{eqnarray*}
so, we have
\begin{eqnarray*}
&&a_1\otimes a_{2(-1)^\prime}a_{3(-1)^\prime}h_{1^\prime}\otimes a_{2(0)^\prime}\otimes a_{3(0)^\prime(-1)}h_{2^\prime1[0]}\otimes a_{3(0)^\prime(0)}h_{2^\prime1[1]}\otimes h_{2^\prime2}\\
&=&a_1\otimes a_{2(-1)^\prime}a_{3(-1)11^\prime}S(a_{3(-1)2})a_{3(0)(-1)^\prime}h_{1^\prime}\otimes a_{2(0)^\prime}\otimes a_{3(-1)12^\prime}h_{2^\prime1[0]}\\
&&\otimes a_{3(0)(0)^\prime}h_{12^\prime[1]}\otimes h_{2^\prime2}.
\end{eqnarray*}

By setting $h=1$ and then applying $\v\o id\o \v\o id\o id\o \v$ to both sides of the above equation, we can prove that Eq.(4.2) holds.
\hfill $\square$

\vspace{3mm}

\noindent{\bf Remark 4.4} (1) Assume that $H$ is a commutative Hopf algebra. Then, by Example 2.3, we know that $(H, \Delta, \Delta^\prime=\Delta^{coH})$ is a commutative Hopf brace.

Suppose that $(A,H,\rho,\varphi)$ is a Hopf matched pair, and the coaction of left $H_{\Delta^\prime}$-comodule bialgebra on $A$ is trivial. Then, Eq.(4.2) and Eq.(4.3) are satisfied. So, according to Proposition 4.3, the bicrossed coproduct $(A\bowtie H, \widetilde{\Delta}, \widehat{\Delta})$ is a Hopf brace, where
$\bar{\Delta}=\widehat{\Delta}$ since the coaction of left $H_{\Delta^\prime}$-comodule on $A$ is trivial.

(2) Let $H$ be a finite dimensional cocommutative Hopf algebra with antipode $S$, ${h_{i}}$ a basis of $H$ and ${h^{*}_{i}}$ the corresponding dual basis of $H^{*}$, and let
$$
R=h_i\o h^{*}_i\in H^{op}\o H^{*}.
$$

Then, by \cite{H. X. Chen}, $R$ is a weak $R$-matrix of $H^{op}\otimes H^{*}$ with the inverse $R^{-1}=S^{-1}(h_i)\o h^{*}_{i}$. So, by Example 3.2 (2), we know that $(H^{op}, H^{*}, \r, \varphi)$ is a Hopf matched pair, and hence one can form the bicrossed coproduct $H^{op}\bowtie H^{*}$, whose comultiplication of $H^{op}\bowtie H^{*}$ is given by
$$
\widetilde{\Delta}(x\o f)=x_1\o h^{*}_{i}f_1h^{*}_{j}\o S^{-1}(h_{j})x_2h_{i}\o f_2
$$
for all $x\in H^{op}, f\in H^{*}$.

 Therefore, according to (1), the dual $(D(H)^{*}, \widetilde{\Delta}, \widehat{\Delta})$ of Drinfel'd double $D(H)$ is a Hopf brace.

(3) Let $H_4=k\{1,g,x,gx\}$ be Sweedler's 4-Hopf algebra with char$k\neq$ 2. As an algebra, $H$ is generated by $g$ and $x$ with relations
$$
g^2=1, ~~x^2=0, ~~xg=-gx.
$$
The coalgebra structure and antipode are determined by
\begin{eqnarray*}
&\D(g)=g\o g, ~~\D(x)=x\o g+1\o x,\\
&\v(g)=1,~~\v(x)=0, ~~S(g)=g,~~S(x)=gx.
\end{eqnarray*}

Let $A=kZ_2$, where $Z_2$ is written multiplicatively as $\{1,a\}$, and
$$R=\frac{1}{2}(1\o 1+1\o a+g\o 1-g\o a)\in H\o A.
$$

Then, by \cite{H. X. Chen}, one can easily see that $R$ is a weak $R$-matrix of $H\otimes A$ with $R^{-1}=R$. So, by Lemma 1.3 in \cite{H. X. Chen}, we have the bicrossed coproduct $H_4\bowtie kZ_2$, and according to (1), we know that $(H_4\bowtie kZ_2, \widetilde{\Delta}, \widehat{\Delta})$ is a Hopf brace.

 (4) Let $A$ and $H$ be two commutative Hopf algebras, and $(A,\r')$ a left $H_{\D^{coH}}$-comodule bialgebra. Suppose that the coaction of the left $H$-comodule algebra on $A$ is trivial. Then, it is easy to see Eq.(4.2) holds, and Eq.(4.3) amounts to
$$
h_{1[0]}\otimes h_2\otimes h_{1[1]}=h_{1[0]}\otimes h_{1[1](-1)^\prime}h_2\otimes h_{1[1](0)^\prime}\eqno(4.4)
$$
for all $h\in H$.

And by Definition 3.1, $(A,H,\rho,\varphi)$ is a Hopf matched pair if $(H,\varphi)$ a right $A$-comodule bialgebra. It is obvious that $(H, \D, \Delta^\prime=\D^{coH})$ is a Hopf brace. So, according to Proposition 4.3, $(A\bowtie H, \widetilde{\Delta}, \bar{\Delta})$ is Hopf brace, whose comultiplication $\widetilde{\Delta}$ is given by
$$
\widetilde{\Delta}(a\otimes h)=a_{1}\otimes h_{1[0]}\otimes a_{2}h_{1[1]}\otimes h_2.
$$

Note that the comultiplication $\widetilde{\Delta}$ of the bicrossed coproduct $A\bowtie H$ is actually the comultiplication of the usual smash coproduct on $A\otimes H$.

\vspace{6mm}

\makeatletter
\newcommand{\adjustmybblparameters}{\setlength{\itemsep}{0\baselineskip}\setlength{\parsep}{0pt}}
\let\ORIGINALlatex@openbib@code=\@openbib@code
\renewcommand{\@openbib@code}{\ORIGINALlatex@openbib@code\adjustmybblparameters}
\makeatother


\begin{thebibliography}{99}
\bibitem{A. L. Agore} A. L. Agore. Constructing Hopf braces. arXiv: 1707.03033, 2017.

\bibitem{A. Smoktunowiczi}  A. Smoktunowicz. On Engel groups, nilpotent groups, rings, braces and the Yang-Baxter
equation. Trans. Amer. Math. Soc. 370: 6535-6564, 2018.

\bibitem{A. Smoktunowiczi1} A. Smoktunowicz and L. Vendramin. On skew braces (with an appendix by N. Byott and L.
Vendramin). J. Comb. Algebra, 2(1): 47-86, 2018.

\bibitem{A. Cap} A. Cap, H. Schichl, J. Vanzura. On twisted tensor product of algebras, Comm. Algebra, 23: 4701 ¨C 4735, 1995.


\bibitem{D. Bachiller} D. Bachiller. Classification of braces of order $p^3$ . J. Pure Appl. Algebra, 219(8): 3568-3603,
2015.

\bibitem{D. Bachiller2016} D. Bachiller. Counterexample to a conjecture about braces. J. Algebra, 453: 160-176, 2016.

\bibitem{D. Bachiller2018} D. Bachiller. Extensions, matched products, and simple braces. J. Pure Appl. Algebra,
222(7): 1670-1691, 2018.

\bibitem{D. Bachiller1} D. Bachiller, F. Ced\'{o}, and E. Jespers. Solutions of the Yang-Baxter equation associated with
a left brace. J. Algebra, 463: 80-102, 2016.

\bibitem{D. Bachiller2}D. Bachiller, F. Ced\'{o}, E. Jespers, and J. Okninski. Iterated matched products of finite braces
and simplicity; new solutions of the Yang-Baxter equation. Trans. Amer. Math. Soc. 370: 4881-4907, 2018.

\bibitem{D. Bachiller3} D. Bachiller, F. Ced\'{o}, and L. Vendramin. A characterization of finite multipermutation solutions of the Yang-Baxter equation. Accepted for publication in Publ. Math., arXiv: 1701.09109, 2017.

\bibitem{F. Ced}  F. Ced\'{o}, E. Jespers, and J. Okni\'{n}ski. Braces and the Yang-Baxter equation. Comm. Math.
Phys., 327(1): 101-116, 2014.

\bibitem{H. X. Chen} H. X. Chen. Quasitriangular structures of bicrossed coproducts. J. Algebra, 204(2):504-531, 1998.


\bibitem{I. Angiono} I. Angiono, C. Galindo, L. Vendramin. Hopf braces and Yang-Baxter operators. Proc. Amer. Math.
Soc. 145(5): 1981-1995, 2017.

\bibitem{J. L¨¹}J. L\"{u}, X. Wang, G. Zhuang. Universal enveloping algebras of Poisson Hopf algebras. J. Algebra, 426:92-136, 2015.

\bibitem{L. Guarnieri} L. Guarnieri and L. Vendramin. Skew braces and the Yang-Baxter equation. Math. Comp. 86(307): 2519¨C2534, 2017.

 \bibitem{S. Caenepeel} S. Caenepeel, B. Ion, G. Militaru, S. Zhu. The factorization problem and the smash biproduct of algebras and coalgebras. Algebr. Represent. Theory, 3: 19¨C42, 2000.

 \bibitem{S. Majid} S. Majid. Physics for algebraists: non-commutative and non-cocommutative Hopf algebras by a bicrossproduct construction. J. Algebra, 130: 17¨C64, 1990.

\bibitem{S. Majid1} S. Majid. Matched pairs of Lie groups and Hopf algebra bicrossproducts, Nuclear Physics B 6: 422¨C424, 1989.

    \bibitem{Sweedler}M. E. Sweedler. Hopf algebras. New York: Benjamin, 1969.

    \bibitem{S. Montgomery} S. Montgomery. Hopf algebras and their actions on
rings. CBMS 82, Providence, RI, AMS, 1993.

\bibitem{T. Gateva-Ivanova}  T. Gateva-Ivanova. Set-theoretic solutions of the Yang-Baxter equation, braces, and symmetric groups.  J. Symb. Comput. 42(11): 1079-1112, 2015.\bibitem{Nakajima1975} A.Nakajima, On generalized Harrison cohomology and Galois object, {\it Math. J. Okayama Univ.}, {\bf 17} (1975):135-148.

\bibitem{T. Br} T. Brzezi\'{n}ski. Trusses: between braces and rings. arXiv: 1710.02870, 2017.

\bibitem{W. Rump1} W. Rump. Modules over braces. Alg. Disc. Math., 2: 127-137, 2006.

 \bibitem{W. Rump2} W. Rump. Braces, radical rings, and the quantum Yang-Baxter equation. J. Algebra,
307(1): 153-170, 2007.

\bibitem{W. Rump2007}W. Rump. Classification of cyclic braces. J. Pure Appl. Algebra, 209(3): 671-685, 2007.

 \bibitem{W. Rump2014} W. Rump. The brace of a classical group. Math. Note, 34(1): 115-144, 2014.

\bibitem{Y. Chen} Y. Y. Chen, Z. W. Wang and L. Y. Zhang. FS-coalgebras and crossed coproducts. Georgian Math. J., DOI: 10.1515/gmj-2017-0037.

\bibitem{Y. Kosmann-Schwarzbach} Y. Kosmann-Schwarzbach, F. Magri. Poisson-Lie groups and complete integrability. I. Drinfel¡¯d bialgebras, dual extensions and their canonical representations, Ann. Inst. H. Poincar Phys. Thor., 49: 433¨C460, 1988.

\bibitem{Zhang2006}  L.Y. Zhang, Long bialgebras, dimodule algebras and
quantum Yang-Baxter modules over Long bialgebras, Acta Math.
Sin. (English Series), 22: 1261-1270, 2006.




\end{thebibliography}
\end{document}